# A Numerical Study of Heat Source Reconstruction for the Advection-Diffusion Operator: A Conjugate Gradient Method stabilized with SVD

**Jing YE, Laurent FARGE, Stéphane ANDRE, Alain NEVEU[†]**

Laboratoire d'Energétique et de Mécanique Théorique et Appliquée (LEMTA),
UMR 7563 CNRS
2 avenue de la Forêt de Haye ,TSA 60604-54518 Vandoeuvre-Lès-Nancy, France

[†] Laboratoire de Mécanique et Energétique d'Evry (LMEE), Université d'Evry Val d'Essonne
40 rue du Pelvoux CE1455, Courcouronnes, 91020 Evry Cedex, France

Corresponding authors:
e-mail: jing.ye@univ-lorraine.fr, stephane.andre@univ-lorraine.fr

**Abstract**

In order to better understand micromechanical phenomena such as viscoelasticity and plasticity, the thermomechanical viewpoint is of prime importance but requires calorimetric measurements to be performed during a deformation process. Infrared imaging is commonly used to this aim but does not provide direct access to the intrinsic volumetric Thermomechanical Heat Sources (THS). An inverse method is needed to convert temperature fields in the former quantity. The one proposed here relies on a diffusion-advection heat transfer model. Advection is generally not considered in such problems but due to plastic instabilities, a heterogeneous and non-negligible velocity field can play a role in the local heat transfer balance. Discretization of the governing equation is made through appropriate spectral approach. Spatial regularization is then achieved through regular modal truncation. The objective of the inversion process lies in a proper identification of the decomposition coefficients (states) which minimize the residuals. When a Conjugate Gradient Method (CGM) is applied to this nonlinear least square optimization, the use of Karhunen-Loève Decomposition (KLD) or Singular Value Decomposition (SVD) on gradient vectors is shown to produce very good temporal regularization. Two test-cases were explored for noisy data which show that this algorithm performs very well when compared to the Tikhonov penalized conjugate gradient method.

# 1   Introduction

Source reconstruction is a classical topic in inverse methods which is shared in common by many scientific fields offering diverse practical subjects of interest:

- ElectroEncephaloGraphy (EEG) and MagnetoEncephaloGraphy (MEG) in medical (imaging) science, for brain current sources reconstruction due to bioelectrical activity, taking into account fields propagation through head tissues;
- Tsunami source reconstruction in geology, from various seismic data recorded from passing waves (the earthquake location problem);
- Acoustic sources reconstruction to improve noise reduction in aeronautic or automotive industries;
- Pollutant sources detection, in environmental research…

All these reverse approaches differ by their objectives (source localization, source strength estimation, characterization of the transport media…), their geometries (point, line, volumes…), the nature of the mathematical forward problem (diffusion-dispersion-reaction-advection-propagation…), the strategies used to help deal with the inevitable ill-conditioning of the inverse source problem (Isakov, 1990). We discuss in this paper the case of an Inverse Heat Source Problem (IHSP) originating from material science (Auffray et al., 2013; Chrysochoos and Belmahjoub, 1992). Thermomechanical Heat Sources (THS) are produced when submitting a specimen of given material to mechanical testing. The deformation processes occurring at various scales of the microstructure are responsible of thermodynamical effects, reversible or not, which can be monitored through temperature/heat flux measurements. Of course, identifying THS during mechanical tests can help the scientist to validate the thermodynamic framework of a macroscopic model (Andre et al., 2012; Vincent, 2008) or to understand more clearly the micromechanisms of deformation and to predict damaging (Doudard et al., 2010; Boulanger et al., 2004).

In previous works (Renault et al., 2008, 2010), the authors perform THS reconstruction with two different strategies, assuming a pure diffusion of heat transfer. The source term $f(r,t)$ was reconstructed from the following parabolic operator $u_{,t} - k\nabla^2 u = f(r,t)$ both in 2D and 1D geometry. Many other papers can be found dealing with the IHSP for the parabolic heat Partial Differential Equation and each of them proceeds in a different manner. In Auffray et al. (2013) for example, the reciprocity gap concept is used to solve the inverse problem. A specification function method is used to identify the number of hidden point sources, their intensity, temporal and spatial locations, holding time. The function specification was also used in a welding application (Rouquette et al., 2007) where the parameters of Gaussian heat sources were determined through a Levenberg-Marquardt approach. The mollification method, a filtering procedure that is appropriate for the regularization of a variety of ill-posed problems has also been proposed by Murio and co-workers (Yi and Murio, 2004) to handle the IHSP. Differentiations of the left-hand side of the above equation are performed on noisy discrete data according to a convolutive process through a Gaussian kernel and give the source term in a rather direct process. The recent works of Delpueyo et al. (2013) can be considered as an application of the mollification method to THS reconstruction. Some works make use of the Boundary Element Method to recover temporal source functions in the 1D case (Farcas and Lesnic, 2004) or line sources in a 2D system with experimental data produced by IR imaging (Le Niliot et al., 2000). Of course many works rely on iterative strategies. In Liu et al. (2015) an original Lie Group Adaptative Method is presented and applied on restoring a stationary space-dependent source function in a 1D transfer problem. Methods based on the Conjugate Gradient Method (CGM) which is used in the present study, were also worked out to reconstruct either temporal or spatial source function in the 1D problem (Hasanov and Pektaş, 2013, Erdem et al., 2013).

But in real applications such as tensile tests performed at high strain rates, a heat transport contribution can originate from advection. Large local velocities can arise from fast imposed solicitations but also because the material undergoes plastic instabilities (strain localization) which produce strong heterogeneities in the velocity distribution. This means that local Peclet numbers greater than one can be obtained. As a result, true THS estimations require

(i) to use Digital Image Correlation (DIC) techniques in order to measure the displacement fields as a function of time, and

(ii) to address the advection-diffusion problem. It means to develop a new algorithm of THS reconstruction, based now on the diffusion-advection operator $u_{,t} + v \cdot u_{,r} - k\nabla^2 u = f(r,t)$, where both variables $u$ and $v$ are produced by measurements.

Only a few inverse studies report some strategy to comply with this problem and generally concern mass transport processes in a fluid flow. In Maalej et al., 2012, the one-dimension linear advection diffusion equation is considered with constant coefficients to recover a transient pollutant point source. This means that the velocity is considered as constant (unidirectional uniform flow) and the inverse problem produces the location and intensity of the source. It is based on an impedance approach i.e. the solution of the transport equation in Laplace domain. Impedances link Laplace concentration at multiple sensor locations with the Laplace transform of the input source intensity. The inverse strategy allows the estimation of the thermal diffusivity and velocity as well. In Rap et al. (2006) or Souza and Roberty (2012), the advection-diffusion-decay equation with constant coefficients and in stationary state is considered for inverse problem reconstruction. After transformation into an Helmhotz equation, inverse approaches are developed around the variational formulation in view of an application where only boundary measurements are available.

Compared with all the cases discussed in this non-exhaustive review, the value of the present work lies in (i) space and time dependent inverse source reconstruction based on transient convection-diffusion equation, (ii) non linear time and space dependent measured velocities used in the inverse process, (iii) no prior information needed on the reconstructed observable (iv) enhancement of the performance of the CGM through a Truncated Singular Value Decomposition (TSVD) approach as stabilizing strategy. As a result the following sketch will be followed:

In section 2, we present the forward problem along with the experimental framework. The heat transfer problem is solved using a spectral method based on a decomposition of the temperature variable on the so-called "branch" basis of orthogonal functions from now on referred to as branch modes. Spectral methods offer a simple way of reducing the model so as to limit ill-conditioning by simply limiting the dimensionality of the basis of functions used. The inversion lies in the identification of the decomposition coefficients (states) through iterative methods. In section 3, an adjoint formulation of the problem is established using an extended Lagrangian in $L^2$ norm for the minimization problem. This leads to a continuous equation which is solved backward to yield directly the driving condition of the Conjugate Gradient Method (CGM). In discretized form, the system will appear highly parameterized and the CGM method is well adapted to handle such cases. In Section 4, an original improvement of the method will be detailed which includes a regularization engine based on a SVD applied to the gradient vectors at each iteration. In Section 5, the performances of the algorithms will be checked through numerical results obtained for two representative test cases.

## 2    Forward Problem

The following 1D transient heat transport problem is considered for the single eulerian space variable $X$ and for time variable $t$.

$$\forall X \in \Omega = ]0;L[ \quad c\frac{\partial T(X,t)}{\partial t} + cv(X,t)\frac{\partial T(X,t)}{\partial X} - k\frac{\partial^2 T(X,t)}{\partial X^2} = q(X,t) \quad (1a)$$

$$\forall t > 0 \quad -k\frac{\partial T}{\partial X}(0,t) = \varphi_1(t) \quad (1b)$$

$$-k\frac{\partial T}{\partial X}(L,t) = \varphi_2(t) \quad (1c)$$

$$T(X,0) = T_0 \quad (1d)$$

Both temperature $T(X,t)$ and local velocity $v(X,t)$ are considered as known (measured) observables. The heat source to be reconstructed is denoted by $q(X,t)$. It depends on both time and space to account respectively for the dynamics of the microstructural evolution of the material and localization effects. Parameters $c$ (heat capacity) and $k$ (thermal conductivity) are supposed constant and perfectly known, because precisely measured in practice for the tested specimen.

In the above system, boundary conditions of Neumann type are considered with "prescribed" heat flux densities $\varphi_1, \varphi_2$. These latter can be either estimated from the temperature gradients at the boundary points $\{0,L\}$ or considered as known. The null heat flux condition for example can be of special interest for the thermomechanical application either because of a symmetry condition on the problem or assuming adiabatic conditions at a sufficient length $L$. The initial temperature $T_0$ (which could also be considered as non-homogeneous over the domain $\Omega$) is also known from measurements of the stabilized temperature field before the experiment starts. According to the real situation, note that convective exchanges occurring at the specimen surface do not appear explicitly here: the 3D→1D dimensional reduction of the problem (fin approximation) combines this heat sink effect to the THS into a single $q(X,t)$ term. The alternative is to consider the advection-diffusion-decay equation directly.

Assuming that $T(X,t)$ belongs to the Hilbert space $L^2(\Omega)$, the energy of the signal remains finite at any time, which means

$$\forall t > 0, \qquad \int_0^L T^2(X,t)dX < \infty \qquad (2)$$

In the frame of a spectral method, the temperature can be decomposed as:

$$T(X,t) = \sum_{i=1}^{\infty} z_i(t)V_i(X) \qquad (3)$$

where $V_i(X)$ are basis functions which depends only on the space variable (principle of separation of time and space independent variables). They are obtained as the solutions of an auxiliary problem (see later). Under the hypothesis that the source $q$ belongs to the Hilbert space $L^2(\Omega)$ as well, it will be decomposed on the same set of modal functions $V_i(X)$ as follows:

$$q(X,t) = \sum_{i=1}^{\infty} b_i(t)V_i(X), \qquad (4)$$

The objective becomes now to reconstruct the source $q$ by identifying the set of decomposition coefficients $B(t) = \{b_i(t)\}$ in a finite space. A truncation of modes is necessary which, as for all spectral methods give a direct entry for regularization (truncation acts as a filter for high frequency modes as will be seen later). In order to introduce a reduced model allowing to calculate the thermal states for a given finite set of modes $N_m$, the modes $V_i(X)$ must be specified. In this work we have considered two distinct sets of basis functions. We chose first to work with the branch basis (Videcoq, et al., 2008) which corresponds to the eigenfunctions of the auxiliary equation of (1-a), along with the generalized Steklov boundary conditions (see equations A-1 and A-2 in appendix A). The theoretical advantage of using such a basis is that, compared to other eigenvalue problems, this one does not consider fixed type of boundary conditions (the eigenvalues appear also in the boundary conditions). As a result, the set of basis functions contains a special family of modes (surface modes) which are specially designed to go along with various boundary conditions. The second set is based on the classical Fourier's modes, adapted to prescribed boundary conditions (see equation B-1 in Appendix B). The use of both sets in parallel in the inverse problem is one feature of this paper so as to see if a gain in efficiency is obtained with the branch set. Irrespective of these choice options, the development of the inverse algorithm follows the exact same path.

Considering the homogeneous linear eigenvalue problem (the velocity for example is considered as constant), the whole system can be rewritten in matrix form with $N_m$ $1^{st}$ order differential equations in terms of state functions (see appendix C):

$$C\dot{Z}(t) = A(t)Z(t) + M(t) + DB(t) = f(Z(t), B(t), t) \qquad (5)$$

$$Z(0) = Z_0 \qquad (6)$$

where capital boldface letters correspond to the vector notation of the n-uplet of numbers denoted by the corresponding small letter as

$$Z(t) = (z_1(t), z_2(t), \ldots, z_{N_m}(t))^T: \quad N_m \times 1$$

$$B(t) = (b_1(t), b_2(t), \ldots, b_{N_m}(t))^T: \quad N_m \times 1$$

The heat capacity matrix $C$ and convection-diffusion matrix $A$ are real-valued with reduced dimension: $\dim A = \dim C = N_m \times N_m$. $C$ is positive and symmetric. $D$ is the matrix related to the scalar product of modes, also of dimension $N_m \times N_m$. $Z_0$ is the projection of the initial temperature $T_0$, obtained from $Z_0 = \langle V, T_0 \rangle_\Omega$, where $V = (V_1(X), V_2(X), \ldots, V_{N_m}(X)): 1 \times N_m$.

In the inverse problem formulation considered in the next section, experimental states $Z(t)$ will be considered as observables. They are the projection of the temperature field data on the $V_i(X)$ modes i.e. $Z(t) = \langle V, T \rangle_\Omega$.

## 3 Inverse Problem Formulation

The least-square criterion used to express the distance between the observable $T$ and the model $T$ is:

$$J(q) = \frac{1}{2}\left(\int_0^{t_f} \|T(M,t) - T(M,t;q)\|^2 \, dt\right) + \frac{1}{2}\|T(M,t_f) - T(M,t_f;q(t_f))\|^2 \tag{7}$$

It corresponds to the Bolza (Loewen et al. 1997) form of the Ordinary Least Square criterion (Hanson, 2006, Liberzon, 2011) as it includes a terminal cost which is generally not considered in thermal inverse problems. According to the choice made in the forward problem (decomposition on modal basis), Equation (7) can be substituted for a criterion expressed in terms of a distance between the set of experimental states components $Z(t)$ and those stemming from the model $Z(t;B)$ carrying the dependency in the unknown state functions $B$. Another consequence of the forward solving method is that the dimensions of vectors $Z(t,B)$ and $B(t)$ are the same, equals to $N_m$, the total number of modes.

We define then the criterion to be optimized as

$$J(Z,B) = \frac{1}{2}\left(\int_0^{t_f}(Z(t) - Z(t;B))^T(Z(t) - Z(t;B))\,dt \right.$$
$$\left. + (Z(t_f) - Z(t_f;B))^T(Z(t_f) - Z(t_f;B))\right) \tag{8}$$

Because the discretized spectral approach of the forward problem leads to a linear system of ODE's (Eq. (5)), this approach leads to a very similar problem as in optimal control (Hanson, 2006). $B$ figures the input and $Z$, the dynamic variable governed by an equation of the type $\dot{Z} = f(Z, B)$ with prescribed initial state $Z(0) = Z_0$.

The goal of the inverse problem is to seek for a feasible trajectory $(Z, B)$ which minimizes J, the cost functional (eq. 8). One difference with respect to optimal control problems is that here, it is aimed this functional be as close to zero as possible because the objective is to reduce a gap between experimental data and model data.

Because solving this problem requires an estimation of the descent gradient $\nabla J$ in the case of a large size problem, the method of Lagrange multipliers or of the adjoint state is used. It is almost always referred to as the Conjugate Gradient Method (CGM). It offers an efficient way of estimating the iterated gradient estimates while considering the dynamical evolution law of the system (eqs. 5,6) as a constraint between the solicitation $B$ and the states $Z$.

The Lagrangian (or Hamiltonian) formulation is constructed as

$$\mathcal{L}(\boldsymbol{Z},\boldsymbol{B},\boldsymbol{\mu}) = J(\boldsymbol{Z},\boldsymbol{B}) + \int_0^{t_f} \boldsymbol{\mu}^T(t)[\boldsymbol{g}(\boldsymbol{Z}(t),\boldsymbol{B}(t),t)]\,dt \tag{9}$$

where $\boldsymbol{\mu}^T(t)$ is the $m$-vector Lagrange multiplier, noted by $\boldsymbol{\mu}(t) = \left(\mu_1(t), \mu_2(t), \ldots, \mu_{N_m}(t)\right)$ also called the adjoint state and where $\boldsymbol{g}(\boldsymbol{Z}(t), \boldsymbol{B}(t), t)$ corresponds to

$$\boldsymbol{g}(\boldsymbol{Z}(t), \boldsymbol{B}(t), t) = C\dot{\boldsymbol{Z}}(t) - \boldsymbol{f}(\boldsymbol{Z}(t), \boldsymbol{B}(t), t) \tag{10}$$

Eq. (10) corresponds to equation (5) modeling the dynamics of the system which in the minimization problem works as an associated constraint: $\min_B \mathcal{L}(\boldsymbol{Z}, \boldsymbol{B}, \boldsymbol{\mu})$ subject to $\boldsymbol{g}(\boldsymbol{Z}, \boldsymbol{B}, t) = \boldsymbol{0}$

The set $(\boldsymbol{Z}, \boldsymbol{B}, \boldsymbol{\mu})$ is considered as independent variables, with in particular

$$\mathcal{L}(\boldsymbol{Z}(t; \boldsymbol{B}), \boldsymbol{B}, \boldsymbol{\mu}) = J(\boldsymbol{B}) \tag{11}$$

Indeed, integrating by parts leads to

$$\begin{aligned}\mathcal{L}(\boldsymbol{B},\boldsymbol{Z},\boldsymbol{\mu}) = \frac{1}{2}&\left(\int_0^{t_f}\big(\boldsymbol{Z}(t)-\boldsymbol{Z}(t;\boldsymbol{B})\big)^T\big(\boldsymbol{Z}(t)-\boldsymbol{Z}(t;\boldsymbol{B})\big)dt\right.\\&\left.+\big(\boldsymbol{Z}(t_f)-\boldsymbol{Z}(t_f;\boldsymbol{B})\big)^T\big(\boldsymbol{Z}(t_f)-\boldsymbol{Z}(t_f;\boldsymbol{B})\big)\right)+\big(C\boldsymbol{\mu}(t_f)\big)^T\boldsymbol{Z}(t_f)\\&-\big(C\boldsymbol{\mu}(0)\big)^T\boldsymbol{Z}(0)-\int_0^{t_f}\left[\big(C\dot{\boldsymbol{\mu}}(t)+A(t)^T\boldsymbol{\mu}\big)^T\boldsymbol{Z}(t)+\boldsymbol{\mu}(t)^TD\boldsymbol{B}(t)\right]dt\end{aligned} \tag{12}$$

For a given set of Lagrange multipliers $\boldsymbol{\mu}$, the differential of the Lagrangian (due to a small change in the control $\boldsymbol{B}$) is

$$d\mathcal{L} = \int_0^{t_f} (-\tilde{\boldsymbol{Z}}(t) + \boldsymbol{Z}(t) - C\dot{\boldsymbol{\mu}}(t) - A(t)^T \boldsymbol{\mu})^T \delta \boldsymbol{Z} dt - \int_0^{t_f} (D^T \boldsymbol{\mu}(t))^T \delta \boldsymbol{B} dt \qquad (13)$$

$$+ \left(C\boldsymbol{\mu}(t_f) - \left(\tilde{\boldsymbol{Z}}(t_f) - \boldsymbol{Z}(t_f)\right)\right)^T \delta \boldsymbol{Z}(t_f) - \left(C\boldsymbol{\mu}(0)\right)^T \delta \boldsymbol{Z}(0)$$

$$= \int_0^{t_f} (-\tilde{\boldsymbol{Z}}(t) + \boldsymbol{Z}(t) - C\dot{\boldsymbol{\mu}}(t) - A(t)^T \boldsymbol{\mu})^T \delta \boldsymbol{Z} dt - \int_0^{t_f} (D^T \boldsymbol{\mu}(t))^T \delta \boldsymbol{B} dt$$

$$+ \left(C\boldsymbol{\mu}(t_f) - \left(\tilde{\boldsymbol{Z}}(t_f) - \boldsymbol{Z}(t_f)\right)\right)^T \delta \boldsymbol{Z}(t_f),$$

since $\boldsymbol{Z}(0) = \boldsymbol{Z}_0$, thus $\delta \boldsymbol{Z}(0) = 0$.

In order to find analytically the gradient vectors, and since we have the freedom to pick $\boldsymbol{\mu}$ in order to make matters simpler, the choice of $\boldsymbol{\mu}$ is made so as to cancel the first integral term and last term of (11) so as to make $d\mathcal{L}$ linear for $\delta \boldsymbol{B}$.

Let $\mu$ be defined through:

$$-C\dot{\boldsymbol{\mu}}(t) = A(t)^T \boldsymbol{\mu}(t) + \tilde{\boldsymbol{Z}}(t) - \boldsymbol{Z}(t) \qquad (14)$$

$$\boldsymbol{\mu}(t_f) = C^{-1} \left(\tilde{\boldsymbol{Z}}(t_f) - \boldsymbol{Z}(t_f)\right) \qquad (15)$$

we have thus the following linear expression

$$d\mathcal{L} = -\int_0^{t_f} (D^T \boldsymbol{\mu}(t))^T \delta \boldsymbol{B} dt = \langle -D^T \boldsymbol{\mu}, \delta \boldsymbol{B} \rangle_{L^2(0,t_f)} \qquad (16)$$

But considering the right term of equation (11), we also have

$$dJ = \int_0^{t_f} \nabla J^T \delta \boldsymbol{B} dt \qquad (17)$$

where $\nabla J$ is the gradient of functional $J$.

Comparing equations (16) and (17) leads to the following gradient equation

$$\frac{d\mathcal{L}}{d\boldsymbol{B}} = \nabla J(\boldsymbol{B}) = -D^T \boldsymbol{\mu}(t) \qquad (18)$$

which is the classical issue of this method. Matrix $D$ is known and corresponds to $\partial \boldsymbol{f}^T/\partial \boldsymbol{B}$. Also classical is the minus sign prefix on the time derivative of the Lagrange multiplier in eq. (14) which requires to solve this ODE backward in time in order to get the $\boldsymbol{\mu}$ vector and as a result, an estimation of the gradient of eq. (18). It is clear that eq.(15) ("initial" condition of the backward problem) would not have been produced without including formally the terminal cost in Eq. (7).

Equations (5),(14),(18) with appropriate initial conditions are called Hamilton's equations and are just the necessary conditions for an interior point optimum of the Lagrangian $\mathcal{L}$ at the optimal set of three vectors $\boldsymbol{Z}^{opt}(t), \boldsymbol{B}^{opt}(t), \boldsymbol{\mu}^{opt}(t)$

We implement the following CGM algorithm as follows, where $\boldsymbol{w}$ will be the notation used for the conjugate gradient vector:

# CGM Algorithm

| Step 1: | Initialization: $n = 0$. |
|---|---|
| | Calculate matrix $A(t), C, D, M(t)$ and initial state vector $\mathbf{Z}_0 = \langle V, T_0 \rangle_\Omega$. |
| | Select initial excitation state vector $\mathbf{B}^{(0)}(t)$. (In practice, $\mathbf{B}^{(0)} = 0$ is chosen to start the algorithm with as less modes as possible). |
| | Calculate first descent vector: $\mathbf{w}^0 = -\nabla J^{(0)}$, where $\nabla J^{(0)} = \nabla J(\mathbf{B}^{(0)})$. |

*Repeat 2 to 7 until one of stopping rules is satisfied while iteration number $n$ augments.*

| Step 2: | Solution of direct problem (5, 6): |
|---|---|
| | initialize $\mathbf{Z}^{(n)}(0) \leftarrow \mathbf{Z}_0, \forall n \in N$ |
| | **for** $k = 0: N_t - 1$ |
| | $C \dfrac{\mathbf{Z}_{k+1}^{(n)} - \mathbf{Z}_k^{(n)}}{\Delta t} = A_{k+1} \mathbf{Z}_{k+1}^{(n)} + M_{k+1} + D \dfrac{\mathbf{B}_{k+1}^{(n)} + \mathbf{B}_k^{(n)}}{2}$      semi-implicit scheme |
| | **end** |
| Step 3: | Solution of adjoint equations: |
| | Initialization of eq. (15) at final time $\boldsymbol{\mu}^{(n)}(t_f) = C^{-1}\left(\mathbf{Z}(t_f) - \mathbf{Z}^{(n)}(t_f)\right)$. |
| | Solution of eq. 14 backward in time |
| | **for** $k = N_t - 1: -1: 0$ |
| | $-C \dfrac{\boldsymbol{\mu}_{k+1}^{(n)} - \boldsymbol{\mu}_k^{(n)}}{\Delta t} = A_k^T \boldsymbol{\mu}_k^{(n)} + \dfrac{(\widetilde{\mathbf{Z}}_{k+1}^{(n)} - \mathbf{Z}_{k+1}^{(n)}) + (\widetilde{\mathbf{Z}}_k^{(n)} - \mathbf{Z}_k^{(n)})}{2}$      semi-implicit scheme |
| | **end** |
| | Calculation of gradient vectors through eq. (18): $\nabla J^{(n)} = \nabla J(\mathbf{B}^{(n)}) = -D^T \boldsymbol{\mu}^{(n)}$ |
| Step 4: | Determine the descent direction: |
| | If $n = 0, \mathbf{w}^n = -\nabla J^{(0)}$; |
| | If $n \geq 1$, calculate $\gamma^n = \dfrac{\|\nabla J^{(n)}\|^2}{\|\nabla J^{(n-1)}\|^2}$ (Fletcher-Reeves criterion), and |
| | $\qquad\qquad\qquad \mathbf{w}^n = -\nabla J^{(n)} + \gamma^n \mathbf{w}^{n-1}$. |
| Step 5: | Line search: |
| | find an optimal step size $\rho^n$, such that $\rho^n = argmin_{\rho \in \mathcal{R}}\{J(\mathbf{B}^{(n-1)} + \rho \mathbf{w}^n)\}$ |
| Step 6: | Calculate next iterate: $\mathbf{B}^{(n+1)} = \mathbf{B}^{(n)} + \rho^n \mathbf{w}^n$ |
| Step 7: | Set $n \leftarrow n + 1$ and return to Step 2 |

The CGM algorithm is proven in the literature to converge to the exact solution within at most $n$ iterations, while $n$ parameters need to be identified. However, noise and/or round-off errors cause the residual to lose accuracy thus leading to lose the orthogonality between the descent vectors and the gradients, especially near the exact solution (Shewchuk, 1994). To overcome this technical issue and to ensure stable convergence, we suggest an original strategy stemming from the computation of conjugate gradients.

## 4   Regularization strategy

*Regularization through the reduction of the model*

One way to constrain the sought for a correct solution is dictated first by the direct approach used here. When spectral methods are used, the best way to achieve a strong spatial regularization effect is by truncation on the number of modes used to reconstruct the solution. Because spectral modes are ordered with increasing spatial frequency, truncation acts exactly like a filter which avoids being sensitive to the measurement noise (generally at high frequency). In other words, eliminating the modes being able to represent the small scale structures of the original temperature profile, makes the model more insensitive to noise. Of course, this tends to produce a smooth quasi-solution of the problem and makes it impossible the recovery of very discontinuous functions. This is quite well known and was largely discussed in our previous works on THS reconstruction in the pure diffusion case (Renault et al., 2010). Many other examples can be found in inverse heat transfer problems where empirical or theoretical (eigenvalue) spectral basis are used to reduce the direct model and regularize the inverse one (Park et al., 2001, Shenefeld et al., 2002). Therefore it will not be discussed further. The way we practically obtain the optimum number of modes (or truncating number $N_m$) requires one of these two conditions to be met:

One criterion is based on the norm of the difference (residuals) between the vector of observed data $T$ and the vector of filtered data $\hat{T}_{N_m}$ which are obtained by a projection of the observed data on the basis of dimension $N_m$. The first criterion is

<u>Criterion 1:</u>   $\tau_{N_m} = \sqrt{\left\| T - T_{N_m} \right\|^2 / (N+1)} < m \qquad m > 1$

It determines the minimum number $N_m$ which makes the "residuals" passing below a threshold level ($\sigma$ = s.t.d of the noise). This complies with the discrepancy principle used generally to stop iterative inversion algorithms.

The second criterion is

<u>Criterion 2:</u> $(\tau_{N_m+1} - \tau_{N_m})/\tau_{N_m} < \varepsilon$

It determines number $N_m$ when the relative decrease of the residuals using the $N_m + 1$ mode is not more than $\varepsilon$.

The truncation of modes allows using a reduced model in the inversion process which is particularly interesting in view of the SVD treatment proposed below.

*Enhancement of the CGM procedure*

In the CGM, regularization is recognized as inherent to the iterative procedure (Hanke, 1995): it is the stopping criterion which definitely avoids further iterating. But because the convergence to the solution is based on the gradient of functional J:

$$\nabla J(t) = \begin{pmatrix} \frac{\partial J}{\partial b_1}(t) \\ \vdots \\ \frac{\partial J}{\partial b_{N_m}}(t) \end{pmatrix} : N_m \times 1 \qquad (20)$$

its approximate estimation can lead to a biased computation of both the descent directions and step size increment, leading to a lack of stability of the inversion and a non optimized convergence rate. Because it obeys to a stochastic character due to the presence of noise, we suggest in this work to implement a Karhunen-Loève Decomposition (KLD) or Truncated Singular Value Decomposition (TSVD) on the Gram matrix of $\nabla J$. Such treatment allows to represent $\nabla J$ with the minimum degree of freedom and then to produce a set of descent directions having the best impact on the optimization process.

Pioneer works which led to the well-known Levenberg-Marquardt method (Pujol, 2007) originate in this idea of trying to adapt the optimal descent vector in both scalar step and descent direction. All these works are based on the explicit calculation of the so-called sensitivity matrix associated to the parameters to be estimated, and generally consider a least-square discretized norm.

For general regression problems of the type $y = F(x): \mathbb{D}(F) \subseteq \mathbb{X} \to F(\mathbb{D}(F)) \subseteq \mathbb{Y}$ where $\mathbb{X}, \mathbb{Y}$ are adequate vector spaces, we can build a priori expectations related to an inverse problem by considering a linearized approach. This is possible if small perturbations on the data yield small perturbations in the model responses. We can consider a linearization of the residuals $R(x) = F(x) - y$ about a solution $x^*$:

$$R(x^* + \Delta x) \approx R(x^*) + S(x^*)\Delta x \tag{21}$$

where $S(x^*)$ is the Jacobian (sensitivity) matrix.

Considering that measurement errors are independent and normally distributed, such that $\tilde{y} = y_{true} + \varepsilon$, where $y_{true} = F(x^*)$ and $E(\varepsilon \cdot \varepsilon^T) = var(\varepsilon)I = \sigma^2 I$, $S(x^*)$ can be used in the case of nonlinear regression to estimate the covariance of the model parameters. Hence, we have

$$Cov(x^*) \approx var(\varepsilon)(S(x^*)^T S(x^*))^{-1} \tag{22}$$

which small magnitude depends on the conditioning of $S(x^*)^T S(x^*)$.

On the other hand, the covariance of the gradient can be estimated from the formulas: $J(x^*) = \|R(x^*)\|^2/2$ and

$$\nabla J(x^*) = S(x^*)^T R(x^*) \tag{23}$$

Thus since $Cov(R(x^*)) = Cov(\varepsilon) = var(\varepsilon)I$, it follows that

$$\begin{aligned}Cov(\nabla J(x^*)) &= E(\nabla J(x^*)\nabla J(x^*)^T) \approx S(x^*)^T E(R(x^*)R(x^*)^T)S(x^*) \\ &= S(x^*)^T E(\varepsilon \cdot \varepsilon^T)S(x^*) = \sigma^2 S(x^*)^T S(x^*)\end{aligned} \tag{24}$$

That means that (22) can be rewritten, without considering the analytical form of $S(x^*)$, as follows:

$$Cov(x^*) \approx \sigma^4 [E(\nabla J(x^*)\nabla J(x^*)^T)]^{-1} \tag{25}°$$

Considering the SVD of $E(\nabla J(x^*)\nabla J(x^*)^T)$, we have

$$E(\nabla J(x^*)\nabla J(x^*)^T) = U(x^*)\Sigma U(x^*)^T \tag{26}$$

where $U(x^*)$ represents the eigenfunction matrix, and $\Sigma$ the diagonal eigenvalue matrix $\Sigma = diag(\lambda_1^2, \lambda_2^2, \ldots, \lambda_n^2)$. As a result the expression for the parameters covariance matrix is:

$$Cov(x^*) \approx \sigma^4 U(x^*)\Sigma^{-1}U(x^*)^T = \sigma^4 \sum_{i=1}^{n} \frac{U_i(x^*)U_i(x^*)^T}{\lambda_i^2} \qquad (27)$$

Since the $\lambda_i^2$ are decreasing, successive terms in this sum make larger and larger contributions to the covariance. If we were to truncate (27), we could actually decrease the variance in our model estimate. Note that because in our case we have $x^*(t)$, the expectancy in the above equations is by definition

$$E\left(\nabla J\left(x^*\right)\nabla J\left(x^*\right)^T\right) = \frac{1}{t_f}\int_0^{t_f} J\left(x^*\right)\nabla J\left(x^*\right)^T dt \qquad (28)$$

However, for our continuous functional (8), using a linearized form allowing an iterative estimation of Gauss-Newton type would require to calculate the sensitivity matrix:

$$S = \left(\frac{\partial z_i(t)}{\partial b_j(t)}\right)_{ij}$$

While the states can be obtained analytically from (5) by making use of Duhamel's formula:

$$Z(t) = \exp\left(\int_0^t C^{-1}A(t')dt'\right)Z_0$$

$$+ \int_0^t \exp\left(\int_{t'}^t C^{-1}A(s)ds\right) C^{-1}(M(t') + DB(t'))\, dt'$$

it is however impossible to calculate matrix S.

It is then preferable to keep on using the gradient vectors $\nabla J(t)$ within the CGM algorithm but to take benefit of a TSVD to improve the search procedure.

Let us define $W$, the energy function associated with $\nabla J(t)$ (Gram matrix) recognized in (28):

$$W = \int_t \nabla J(t)\nabla J^T(t)dt: \quad N_m \times N_m \qquad (29)$$

Because $W$ is a compact and positive (Hilbert-Schmidt) operator, the spectral theorem states that there exists a complete and orthogonal set of eigenvectors associated to $W$ in its Hilbert space $\mathcal{H}$. The real singular values are all positive with $\lambda_1^2 \geq \lambda_2^2 \geq \cdots \geq 0$. Matrix $W$ can be rewritten using SVD:

$$W = U\Sigma U^T = \sum_{i=1}^{N_m} \lambda_i^2 u_i u_i^T \tag{30}$$

with

$$U = \{u_i\}_{i \leq N_m} \qquad \Sigma = diag(\lambda_1^2, \lambda_2^2, \ldots, \lambda_{N_m}^2)$$

with the eigenvectors verifying

$$U^T U = I \tag{31}$$

Any function in the Hilbert space can be projected uniquely on the eigenvectors set of $W$, because $\{u_i\}_{i \leq N_m}$ is a complete family in the definition domain $\mathcal{H}$. We have:

$$\forall t, \forall f \in \mathcal{H}, \qquad f(t) = \sum_{i=1}^{N_m} u_i z_i(t), \tag{32}$$

In particular $\nabla J(t)$ itself is decomposed as

$$\nabla J(t) = \sum_{i=1}^{N_m} u_i z_i(t) \tag{33}$$

with singular values of $\nabla J(t)$ being the eigenvalues of the Gramian matrix $W$. It is easy to see that

$$z_i(t) = u_i^T \nabla J(t) \tag{34}$$

and from (30), (32), we obtain

$$\langle z_i(t), z_j(t) \rangle_{L^2[0,t_f]} = \int_t z_i(t) z_j(t) dt = \delta_{ij} \lambda_j^2 \tag{35}$$

Eq. (33) represents SVD (or KLD) implemented on $\nabla J(t)$. The related term $u_i z_i(t)$ is taken as the $i^{th}$ principal component of $\nabla J(t)$.

As mentioned in (Del Barrio, 2012), the TSVD technique provides an efficient way of capturing the dominant components of a finite dimensional process. Applied to the CGM, it offers a strategy to proceed to a descent operated in those selected directions where the maximum gain will be made (Iterative reduction of the dimensionality of the optimization procedure). Modes of high frequency are removed in order to have consistently behaving solution. A $r$-rank approximation leads to take

$$\nabla J_r(t) = \sum_{i=1}^{r} u_i z_i(t) \tag{36}$$

which offers the best $r$-rank approximation of $\nabla J(t)$ in terms of a least-square norm of the errors as follows:

$$\|e_r\|^2 = \|\nabla J - \nabla J_r\|^2 = \sum_{i=r+1}^{N_m} \sigma_i^2 \tag{37}$$

The rank of truncation $r$ is determined in practice as the threshold under which the first $r$ singular values contribute to nearly 95% of the "energy", i.e. of the total informational content of the gradient matrix

$$\frac{\sum_{i=1}^{r} \sigma_i^2}{\sum_{i=1}^{N_m} \sigma_i^2} \geq 95\%. \tag{38}$$

Consequently, the use of TSVD on $\nabla J(t)$ (eq. 36) at Step 3 of the algorithm allows keeping the useful information in terms of descent direction and makes sure that the filtered gradient vectors are not far from the exact ones. The whole conjugated gradient method stabilized by TSVD leads to a fast convergence. Cabeza et al. (2005) applies a TSVD during the inverse procedure within a sequential future times steps method.

The stopping criterion that was used is

$$J(B^{n+1}) < \varepsilon = \tau \sigma^2 \tag{39}$$

according to the discrepancy principle where $\tau$ is a constant relative to the functional (7).

# 5 Numerical examples, results, discussion

We give first in section 5.1 all information regarding the inversion test cases used to produce the results of inversion presented in section 5.2. These illustrate all the aspects of the CGM algorithm we developed along with the regularization strategies discussed in section 4. We chose to report in this section only those results obtained when using the second set of spectral basis (Fourier – Appendix B) as it is proved at the end of the study (Section 5.3) to be more efficient than the first set (Branch modes – Appendix A).

## 5.1 General features of the simulations of inversion:

In order to state the performances of our inverse approach, two test cases were selected to cover specific behaviors and features of possible THS. The first one (Section 5.1.1) concerns a stationary but spatially dependent heat source, which is not continuously differentiable. The second one (Section 5.1.2) concerns a transient and spatially dependent source, that is continuously differentiable. Different types of boundary conditions have been considered which all lead to the same performances of the method (same quality of the identification process in terms of residuals on both the source and the temperatures). In what follows, anyway, results are reported for Boundary Conditions of the Neuman type only. These are more likely with respect to the experimental conditions we encountered in practice in the thermomechanical tensile test application. In both test cases, a velocity profile is defined (Figure 1, eq. 42), which is more or less in accordance to what was found experimentally (Ye et al., 2015). It depends both on space and time. Note that in real experiments, the velocity function considered as an input in the inverse problem, results from a field measurement which is noisy and potentially biased.

The simulated "experimental" temperature obtained with the model and given input source and velocity field is produced thanks to a Finite Element software (FlexPDE®). Any "inverse crime" is therefore avoided (totally different solving methods were considered for the direct problem and for the inverse one). For the inversion, an artificial white noise $e(X,t)$ has been superposed on the two simulated observables. Figs 2 and 3 give the plots of the simulated temperature profiles which will be used as input in the inversion process for test cases Nr 1 and 2. This noise is considered as additive and obeys a Gaussian distribution characterized by a 0-mean and a constant standard deviation σ. In order to compare inversion results for various levels of noise, the magnitude of the Signal to Noise Ratio (SNR) is calculated for a given time as:

$$SNR = \frac{\max\{T\}}{2\sigma}: \quad \sigma \text{ standard deviation of noise} \tag{40}$$

As the spectral method is intended to separate temporal from spatial variables, two distinct regularizations can be performed for a global regularization of the inverse process.

The temporal regularization is embedded in the proposed modified CGM algorithm which minimizes iteratively a functional defined on the states. Recordings of the measurement noise $e(X,t)$ (in the initial isothermal condition) can be projected onto the set of eigenmodes to get the corresponding perturbation signal on the states $Z^e(t)$. The stopping rule for iterative regularization is chosen according to the following version of the discrepancy principle:

$$J(B) < J_e = \frac{1}{2}\int_0^{t_f} Z^e(t)^T Z^e(t) dt \tag{41}$$

The initial temperature is taken as $\boldsymbol{T}(0) = 0$, which implies $\boldsymbol{Z}(0) = 0$

### 5.1.1. Test Case 1: Stationary Source

This problem originates from a test case found in (Versteeg and Malalasekera, 2007) and used to validate the numerical solution obtained for a direct (forward) problem. The source term is only of class $C^0$ or piecewise $C^1$ on $[0;L]$ with two "particular" points localized at $X = 0.6$ and $X = 0.8$ (unit in cm) as it can be seen on the plots of Figs.(4-5). The following values are considered for the parameters of the model:

| Heat capacity | $c = 1\, J/(cm^3 \cdot K)$ |
|---|---|
| Thermal conductivity | $k = 0.03\, W/(cm \cdot K)$ |
| Length | $L = 1.5\, cm$ |
| Final time | $t_f = 40s$ |
| Boundary fluxes settings | $\varphi_1(t) = -0.005 \exp(0.1742 t)$ $\varphi_2(t) = 0.005 \exp(0.1249 t)$ |

The following mathematical form has been chosen for the velocity field (Figure 1):

$$v(X,t) = 0.1 c \frac{t}{t_f} \tanh(3 \frac{t}{t_f}(X - L/2)) \tag{42}$$

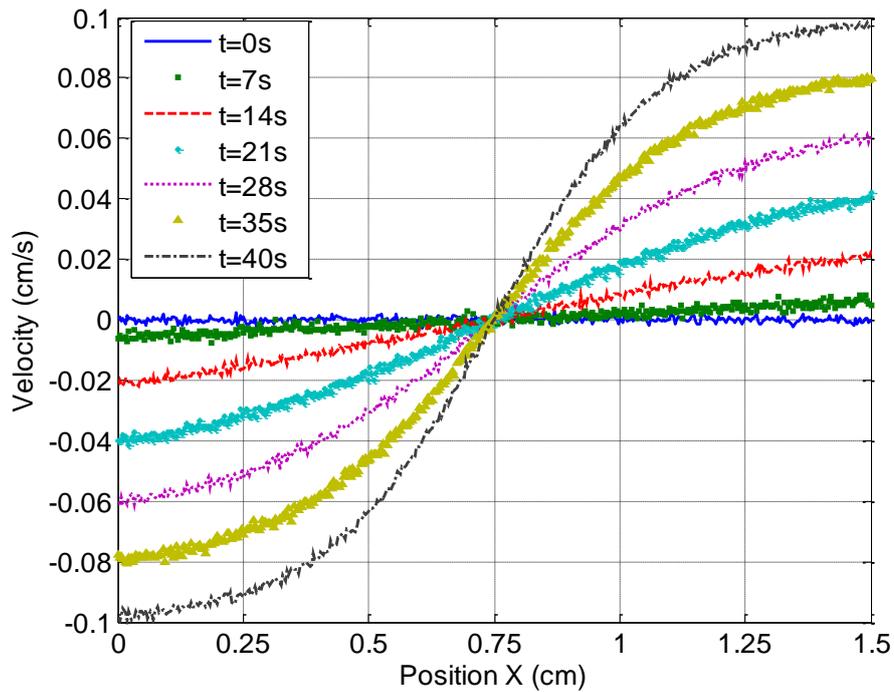

Fig.1. Modeled velocity field (Noise corresponds to 2% of max ($|v|$))

In figure 2, a noise of s.t.d $\sigma = 0.3\, K$ (SNR≈60) is considered for the temperature.

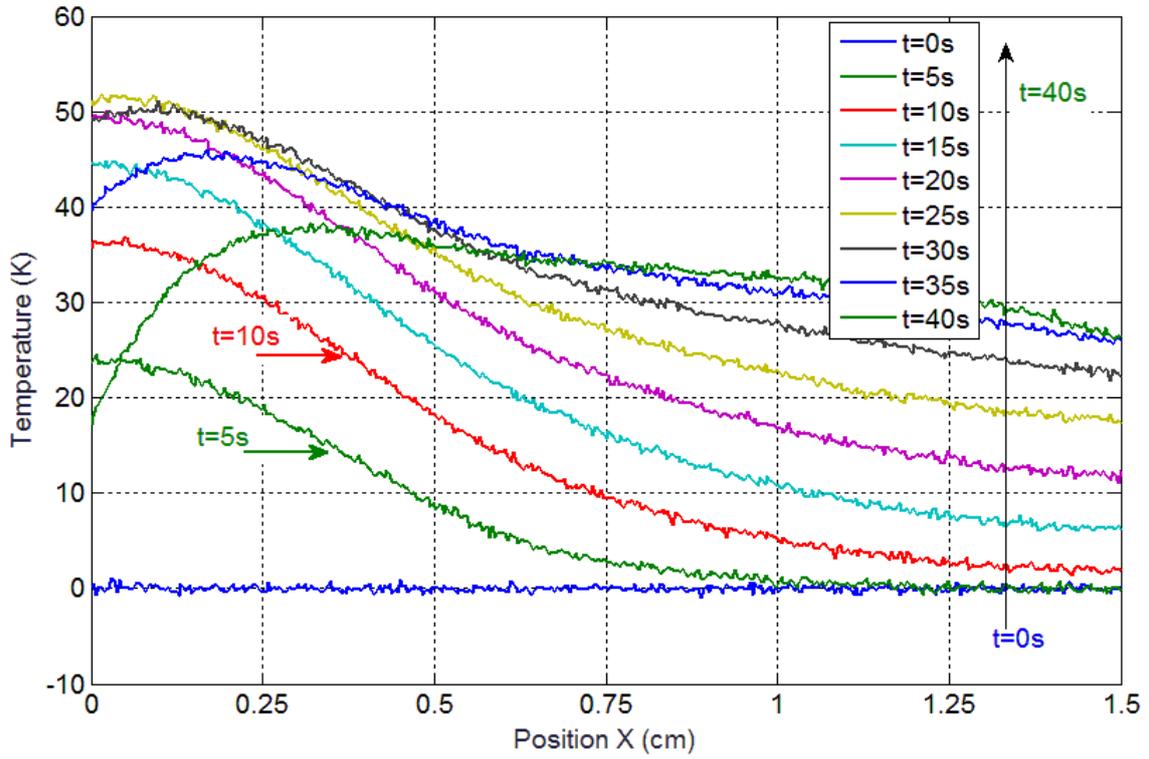

Fig.2. Test-case 1: simulated temperature profiles (Noise $\sigma = 0.3\ K$)

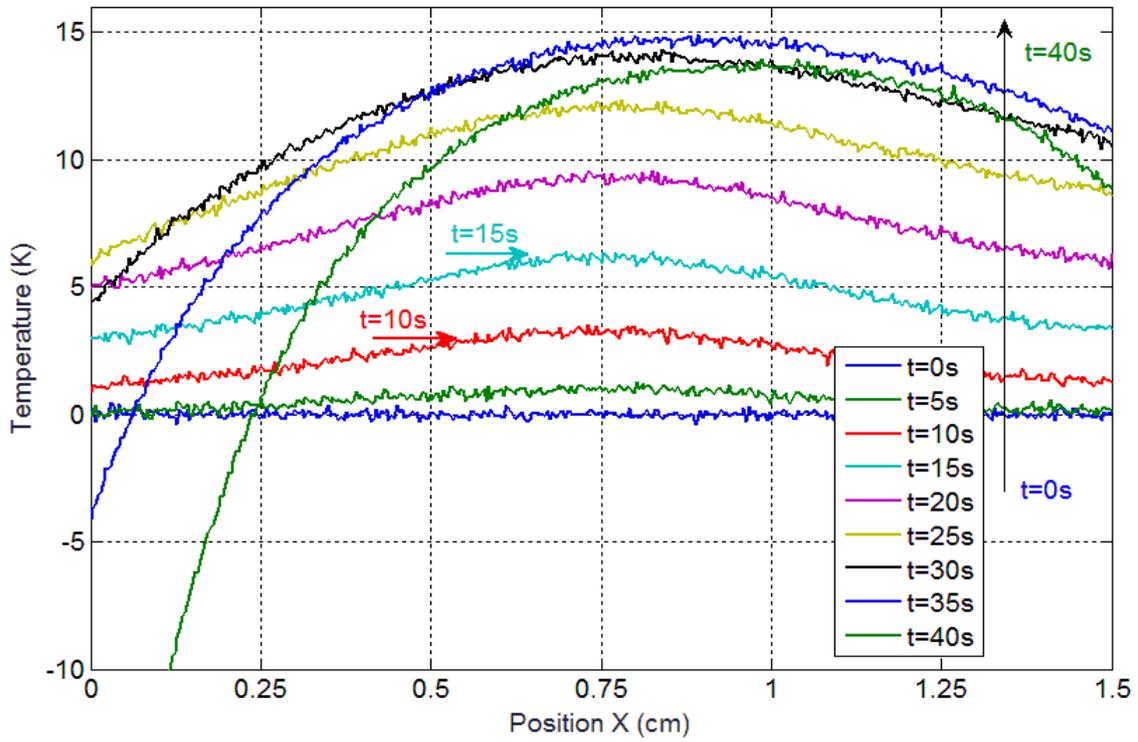

Fig.3. Test-case 2: simulated temperature profiles (Noise $\sigma = 0.13\ K$)

### 5.1.2. Test case 2: Temporal and Spatial varying Source

For the second test-case, a "smooth" but non stationary source is considered

$$q(X,t) = \left[\frac{t}{30} + \sin\frac{t}{10}\right] \exp\left(-\frac{(X-L/2)^2}{0.1}\right) \qquad (43)$$

All other parameters remain unchanged. Figure 3 shows the plots of the transient temperature profiles obtained with an input considered THS given by Eq.(43), the time-dependent velocity profiles shown in Fig 1, and a noise of s.t.d $\sigma = 0.13\ K$ which gives a SNR≈50 (same order of magnitude as for the plots of Fig.2)

### 5.2 Results of inversion:

Before presenting the results, it is of prime importance to have in mind that the two THS sources considered in test-cases 1 (stationary source of Fig. 4) and 2 (varying source starting from a zero value, eq.43) produce very different temperature variations inside the domain for a given time. Because the s.t.d of the noise is determined experimentally by the sensors signals and do not change in general during the experiment, the SNR changes with the time considered for inversion. Therefore, and in view of comparing the behavior of the algorithms we developed, results of inversion are given in general for similar SNR which, in turns, implies to consider either a different s.t.d of the noise or a different time of the reconstruction.

Test-case 1 is first considered. Fig.4 reports the reconstructed THS, the theoretical one (used as input for the forward simulation), and the difference between them, when no noise perturbation is added on the input temperature and velocity data. The thermal source is well rebuilt, thus proving the good working of the algorithm. This reconstruction was obtained for a huge amount of Fourier modes ($N_m = 100$). Since the Fourier spectral method is a global approximation, its main drawback is limited performances of the reconstruction at the limits of the domain or at a singular point, even if more modes are used (an infinity of modes would actually be required).

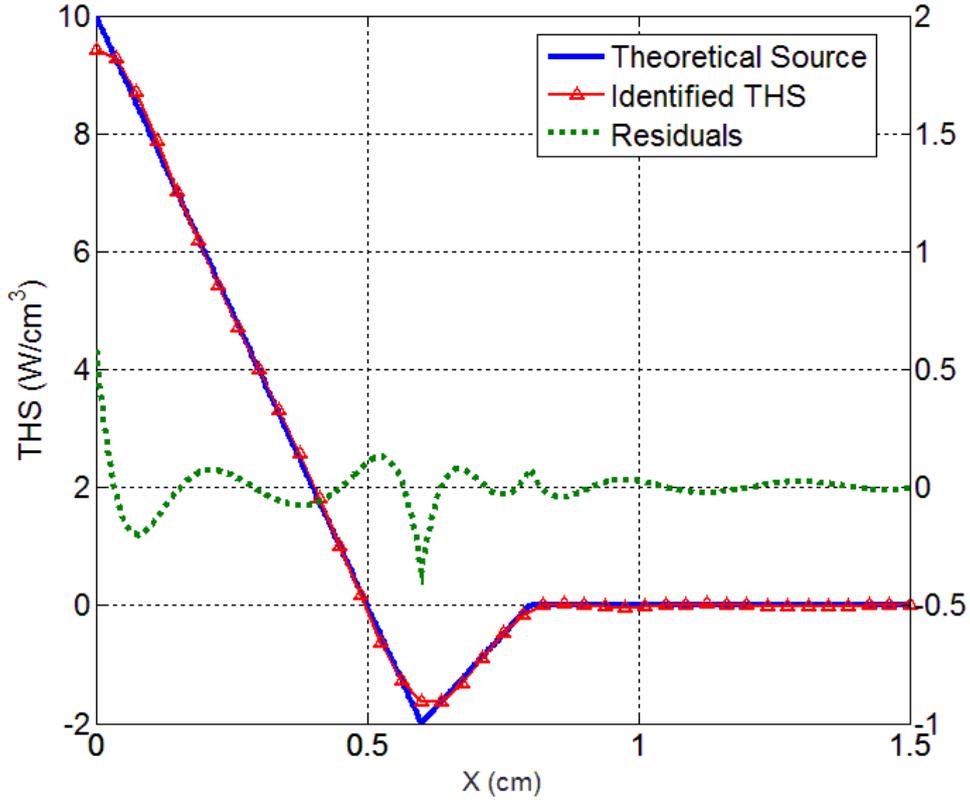

Fig.4. Test Case 1: THS reconstructed at $t = 20\ s$ - No noise - $N_m = 100$ Fourier modes.

Each of the figures 5a-b-c-d plots the same type of results, with on a same column, the results for two different levels of noise ($\sigma = 0.3\ K$ and $\sigma = 0.6\ K$) and on a same line, the results comparing the reconstruction with the regular CGM against those obtained with the present CGM regularized by TSVD. For $\sigma = 0.3\ K$ (SNR$\cong$60), the optimal number of modes is determined according to the algorithm criteria settings and is found to be $N_m = 18$. For $\sigma = 0.6\ K$ (SNR$\cong$30), we keep the same number of modes (equal dimension of spectral basis) to show the efficiency of the TSVD regularization in terms of stabilization of the algorithm (increasing the noise would logically lead to diminish the number of reconstruction modes).

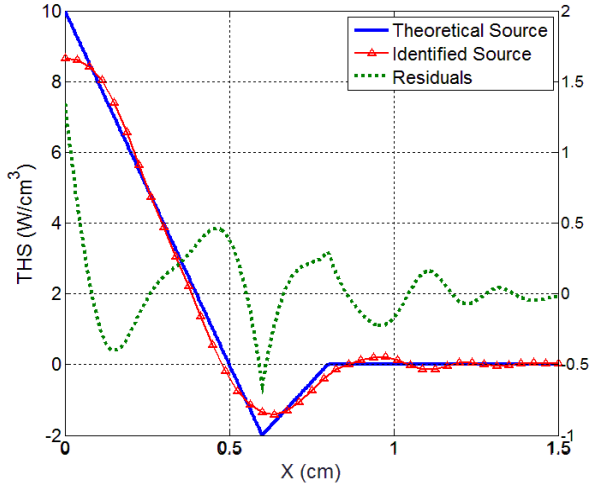
Fig.5a: Noise $\sigma = 0.3\ K$ - without TSVD

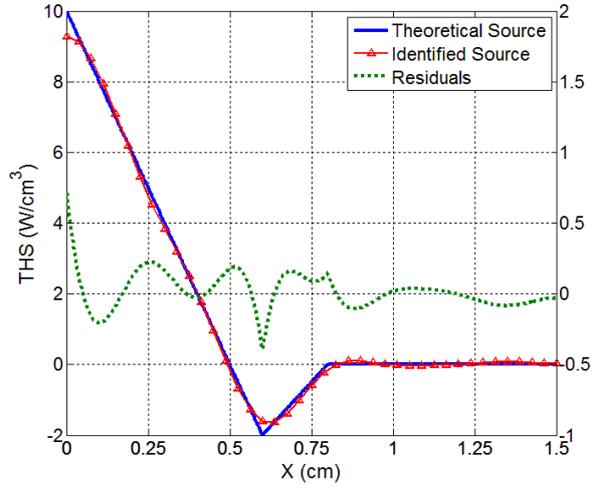
Fig.5b: $\sigma = 0.3\ K$ with TSVD

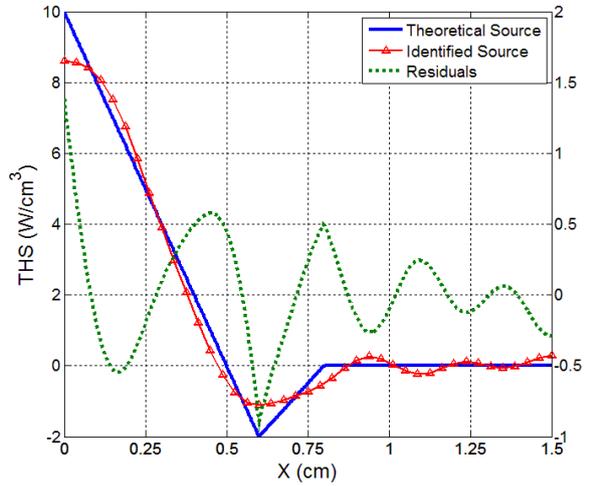
Fig.5c: $\sigma = 0.6\ K$ without TSVD

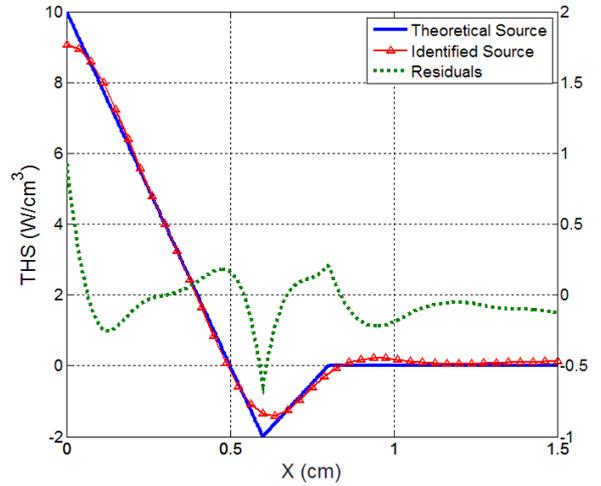
Fig.5d: $\sigma = 0.6\ K$ with TSVD

Fig. 5-abcd. Test Case1: Reconstructed THS at t=20s with $N_m = 18$ modes. Influence of noise and improvement with the TSVD version of the CGM.

The results show clearly that the regular CGM extended with a TSVD approach for optimizing the descent directions produces better results. Apart from what happens at the borders, which is highly constrained by the choice made for the modes (their ability to respect them), the relative errors on the identified THS are diminished by a factor of roughly 2. We consider for this

$$error^{THS} = \frac{\max|q_{id} - q_{exact}|}{\max(q_{exact}) - \min(q_{exact})},$$

which most of the time reflects the discrepancy at the singular point.

When the noise is doubled (SNR decreased by a factor of 2), the regularized CGM (Fig5d) works with a smaller error than the regular CGM with the initial noise (Fig.5a). Figure 6 shows that the reconstruction algorithm performs equally even for different SNR (for the same noise s.t.d, the SNR increases with the experimental time). In all cases, the algorithm produces a very good optimization on the temperature profiles (Fig. 7abcd). The residuals appear non biased with zero mean ($\mu_{res} \ll \sigma_{res}$) where $\sigma_{res}$ is the standard deviation on the residuals, always found roughly equal to the std introduced on the input noise.

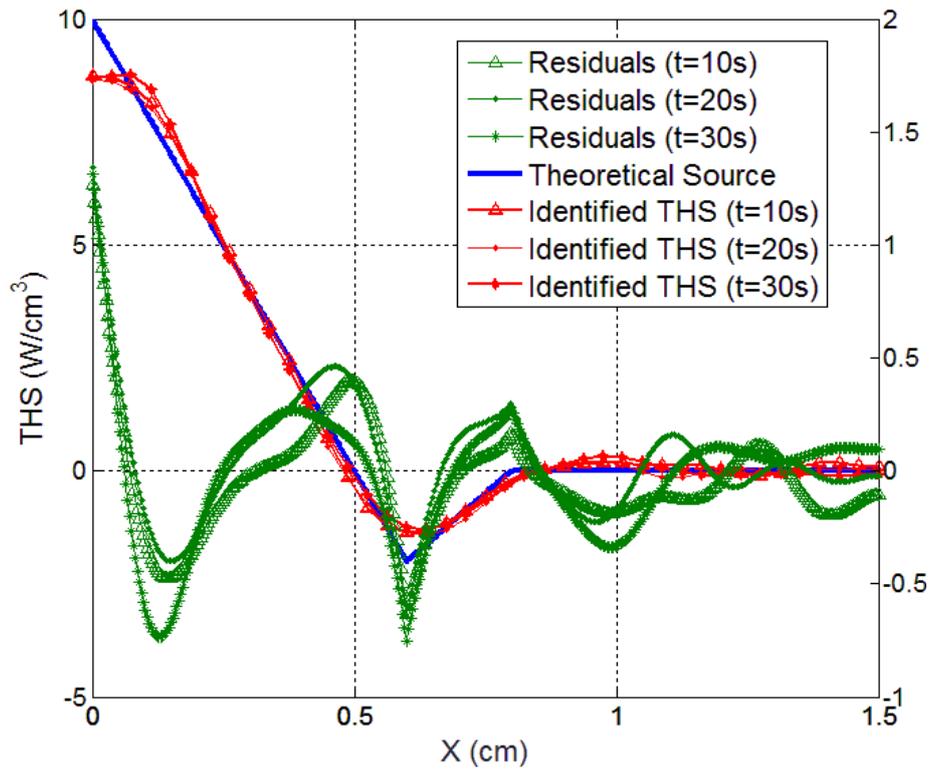

Fig. 6. Test Case1: Reconstructed THS for different experimental times (t=10s, 20s, 30s) with $N_m = 18$ modes, $\sigma = 0.3\ K$ (algorithm without TSVD)

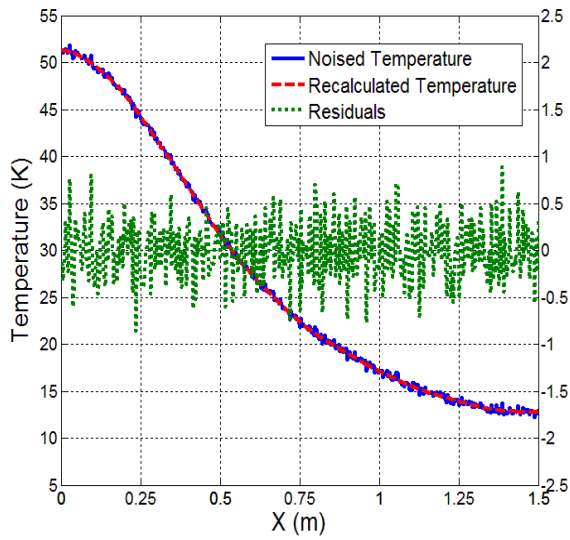 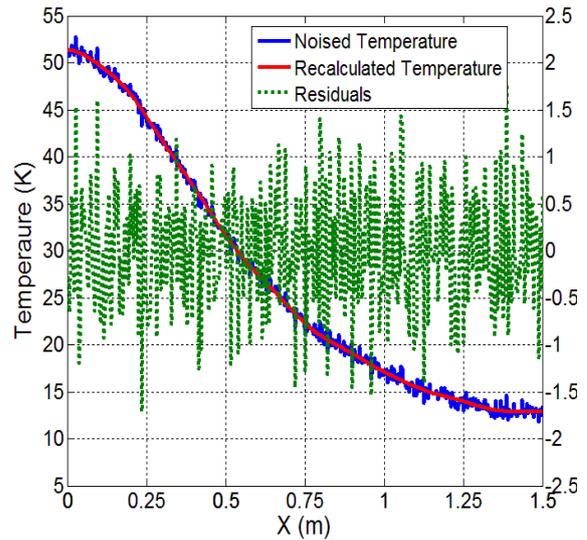

a) $\sigma = 0.3K$ ($\mu_{res} = 0.0143$; $\sigma_{res} = 0.303K$)   b) $\sigma = 0.6K$ ($\mu_{res} = 0.0407$; $\sigma_{id} = 0.605K$)

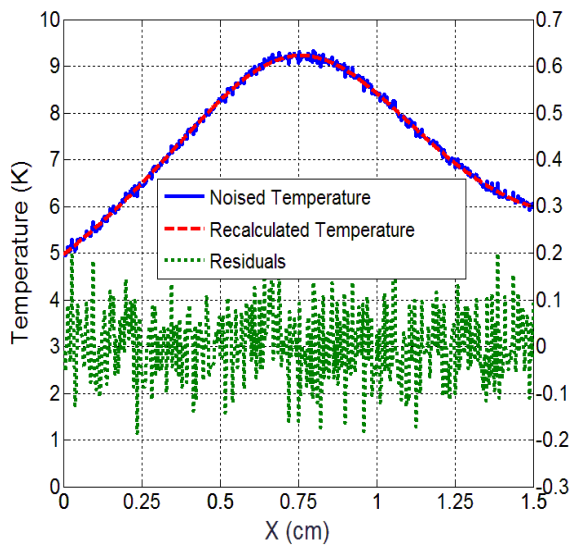 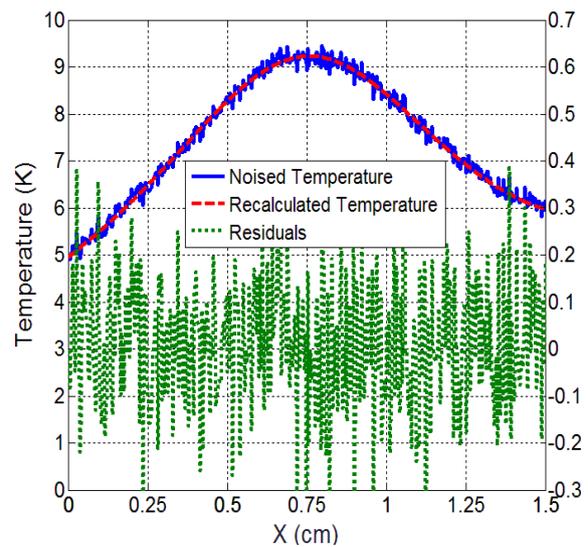

c) $\sigma = 0.07K$ ($\mu_{res} = 0.0047$; $\sigma_{id} = 7.08\,10^{-2}K$)   d) $\sigma = 0.13K$ ($\mu_{res} = 0.0125$; $\sigma_{id} = 0.131K$)

Fig. 7. Input noisy and output reconstructed temperature profiles (left Y axis) – Post-identification residuals (right Y-axis) : (a,b) Test Case 1 - (c,d) Test Case 2. (The mean and s.t.d of the residuals are given in brackets).

Regarding the way both versions of the CGM algorithm work, Figure 8 plots the behavior of functionnal J against the number of the iterations run by the algorithm. It can be seen that the stopping criterion (eq. 41) is reached earlier with the regular CGM (96 iterations) than with the stabilized CGM-TSVD (161 iterations). This indicates that in the regular version, a local minimum is found too early and that the algorithm is trapped there. One very important point to note is that although the maximum number of iterations is augmented, both versions of the CGM algorithm produce the same CPU time. This is because the regularizing effect of TSVD allows less modes to be used for the inversion model which compensates approximately the increase in iterations.

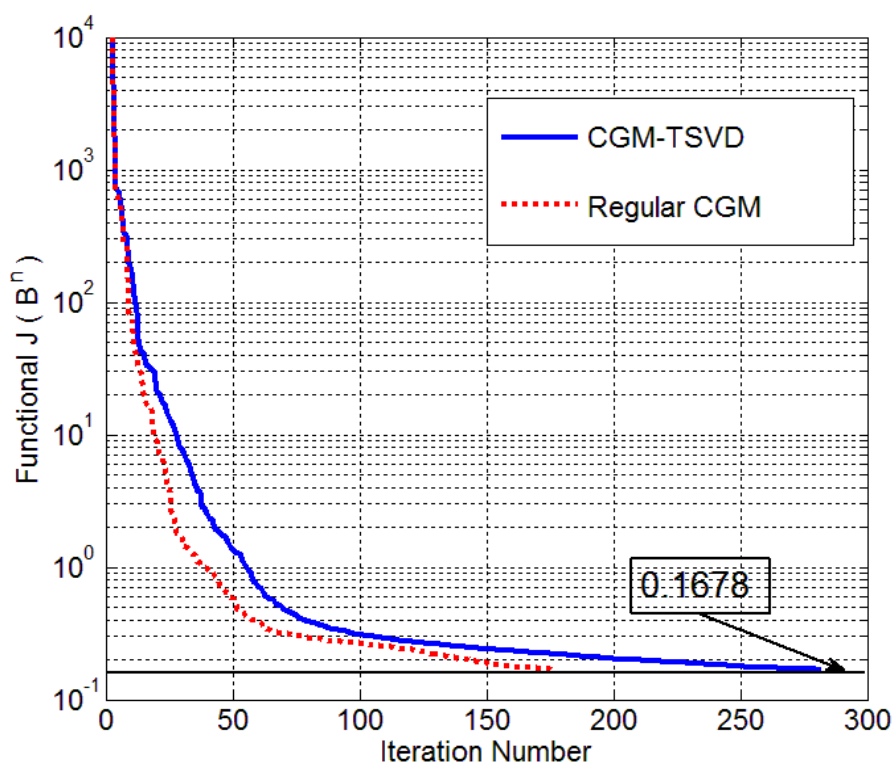

Fig.8. Euclidean norm of the residuals using regular CGM and CGM-TSVD.

We have also compared our results when use is made of a Tikhonov-type regularization applied on CGM (Hochstenbach and Reichel, 2010). We then consider the modified functional

$$J_\alpha(B) = J(B) + \frac{1}{2}\alpha \|B\|^2$$

$$= \frac{1}{2}\left(\int_0^{t_f} \left(\tilde{Z}(t) - Z(t;B)\right)^T \left(\tilde{Z}(t) - Z(t;B)\right) dt \right.$$

$$+ \left(\tilde{Z}(t_f) - Z(t_f;B)\right)^T \left(\tilde{Z}(t_f) - Z(t_f;B)\right)$$

$$\left. + \alpha \int_0^{t_f} (B(t))^T (B(t)) dt \right) \quad (43)$$

where Tikhonov coefficient was set to $\alpha = 0.0001$ (value ensuring empirical systematic convergence). In table 1 are reported the main features of the optimization process which allow comparing the 3 options. The same two different levels of noise are considered. Within the robust criteria used for determining convergence, it can be seen that the proposed technique performs always better than the others. Moreover, one non negligible practical advantage is that the number of truncated singular values is determined in a more straightforward way that the Tikhonov penalization coefficient is. An inappropriate $\alpha$ value could yield a biased solution. According to all these convergent observations, the power of TSVD applied to the descent vectors matrix apparently relies in that it produces a global optimum.

Table 1: Output Parameters of the inversion process at convergence (Test-case 1).

| Noise ampl (K) | Noise s.t.d $\sigma = 0.3K$ | | | Noise s.t.d $\sigma = 0.6K$ | | |
|---|---|---|---|---|---|---|
| | $J_{min}$ | $I_{max}$ | Error THS | $J_{min}$ | $I_{max}$ | Error THS |
| Original CGM | 0.1181 | 96 | 10.1% | 0.4678 | 65 | 14.2% |
| Tikhonov CGM | 0.1191 | 150 | 8.3% | 0.4671 | 82 | 12.5% |
| CGM-TSVD- | 0.1172 | 161 | 5.8% | 0.4778 | 95 | 8.3% |

We now consider test-case 2 which is more likely to resemble that of the application. It can be seen on Figs 9-ab that the results are perfect when CGM-TSVD is used ($error^{THS}$ less than 1.26%). Non perfect reconstruction can be observed on Fig.9-a with regular CGM mainly at the borders and where THS exhibits the higher gradients ($X \approx 0.5 - 1$). Here the noise considered on the input temperature field was of s.t.d. $\sigma = 0.13\ K$ (SNR≈50 at $t = 40s$).

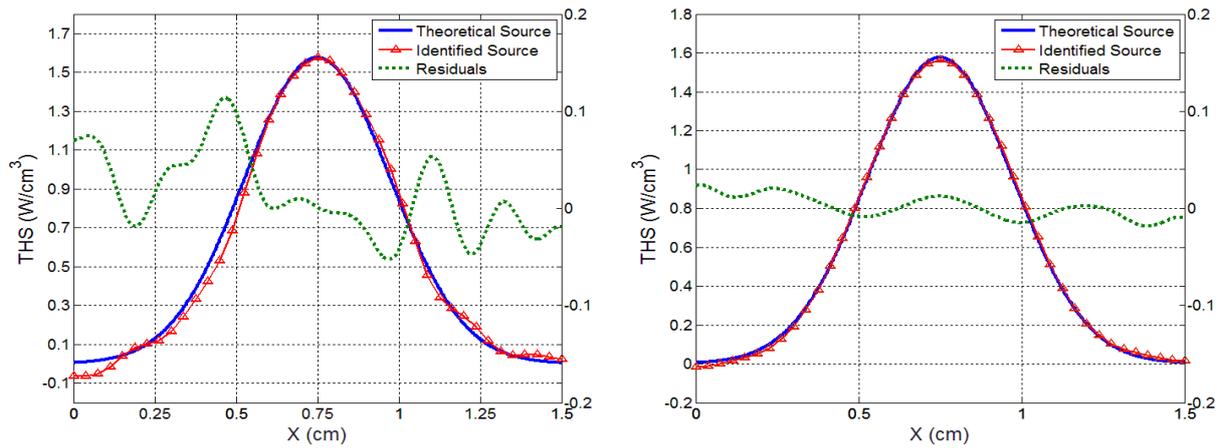

a) without TSVD - $I_{max} = 45$      b) with TSVD - $I_{max} = 44$

Fig. 9. Test case 2 : Reconstructed THS at $t = 20s$ - $N_m = 18$ - $J_{min} = 0.33$

## 5.3 Efficiency of different spectral basis

We now compare the results obtained when use is made of the Branch basis in comparison with the Fourier basis. This will be made only for test-case 1 (stationary source).

We first show in Fig. 10 the 4 first modes corresponding to each basis set. This major difference can be noticed: all Fourier modes have a zero tangent at the limits which is not the case for diffusive Branch modes due to the Steklov condition.

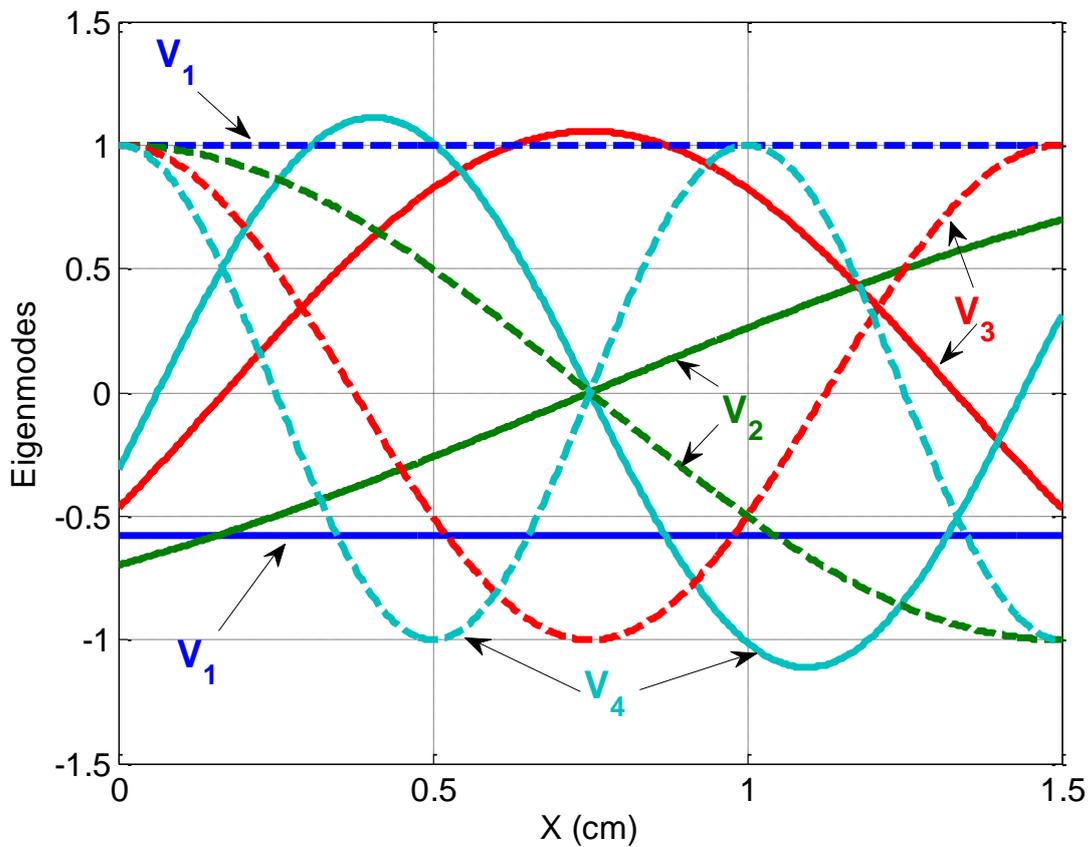

Fig 10. Comparison of the first 4 Branch modes (solid lines) and Fourier modes (dashed lines)

For both choices of modal basis, the initial value of the parameter vector (states $\boldsymbol{B}^{(0)}$) is set to be zero. The results shown columnwise on figures 11-abcdef correspond to 3 levels of noise (un-noisy $= 0$, $\sigma = 0.02$, $\sigma = 0.1$).

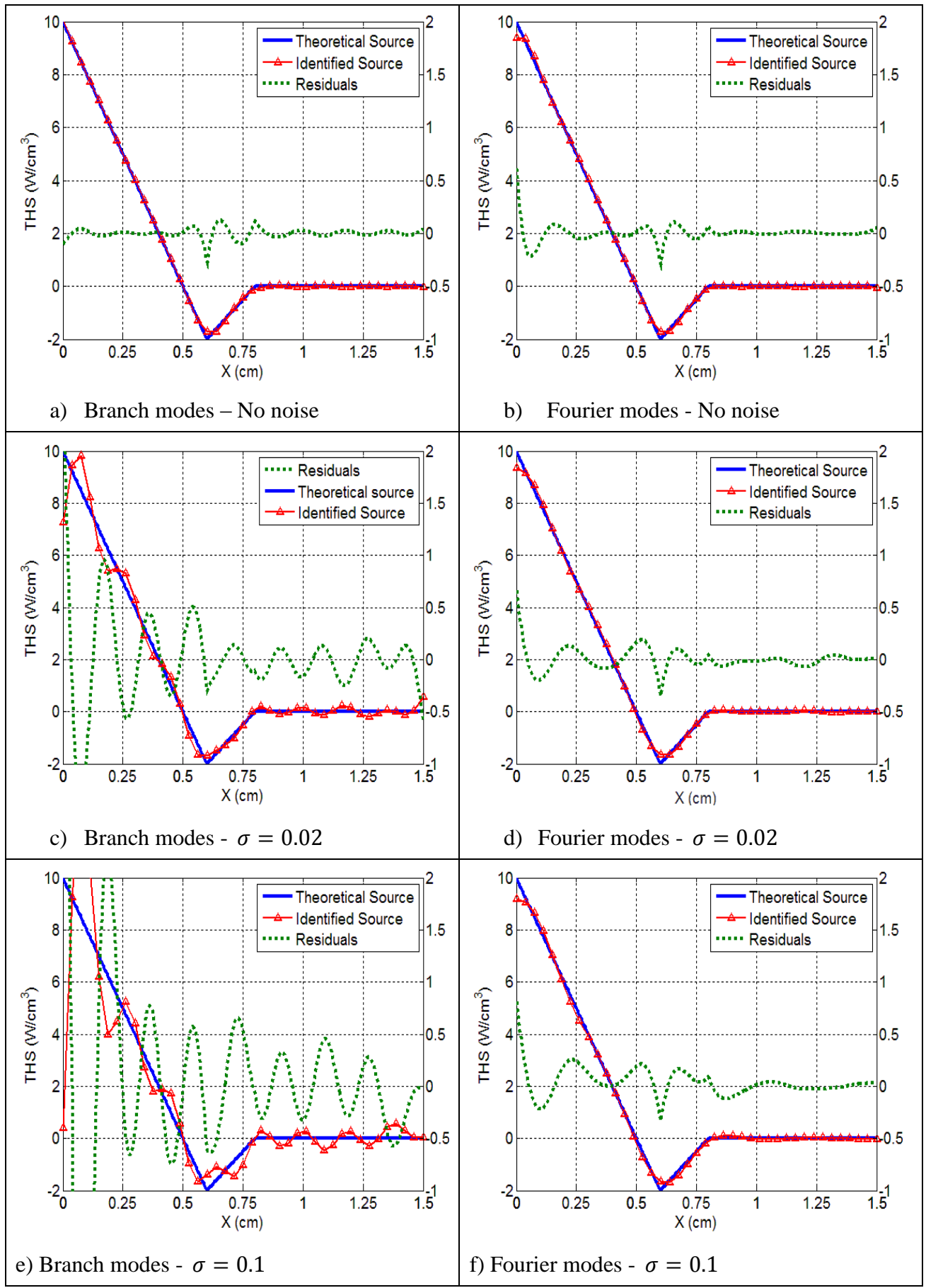

a) Branch modes – No noise
b) Fourier modes - No noise
c) Branch modes - $\sigma = 0.02$
d) Fourier modes - $\sigma = 0.02$
e) Branch modes - $\sigma = 0.1$
f) Fourier modes - $\sigma = 0.1$

Figs. 11-abcdef. Test case 1 - THS reconstructed at $t = 20s$ with 18 modes for both basis

For the Fourier basis, the number of minimum (optimal) modes is $N_m = 18$ for the three levels of noise. The same number of branch modes is considered for the reconstruction of figs 11-ace in order to provide a good basis for comparison because, as explained below, the minimum (optimal) number of branch modes to get good temperature residuals is generally found greater than 18, which has the consequence of introducing instabilities in the reconstruction.

In the absence of noise and for equal number of modes, the Branch basis appears much well-behaved than the Fourier basis especially regarding all crucial locations as the singularity and the frontiers (Figs. 11-ab). It is important to mention here that without any added artificial noise on the simulated temperature, the optimal number of branch modes is found equal to 25. Figure 12 shows clearly that this number is too large with respect to the quality of the reconstruction. While the general trend is correct, the reconstructed signal is covered by oscillations. These are due obviously to the high frequency modes for which the associated states are found to be very sensitive to any small perturbations on the input signal, as we are considering here the case of a pure computational noise resulting for the direct model. This has been proved as related to the way the spectral basis concentrates the major part of the power spectral density of the signal to be reconstructed in the first few modes (Barrio et al., 2012). Comparing Fig.12 to Fig. 11-a allows to understand this. Reducing the number of branch modes allows reaching a perfect source reconstruction, without altering here the residuals on the temperature. In figures (11-cd) and (11-ef) we can see that as soon as a small amount of noise is introduced in the data, the Branch basis method becomes much less efficient and can even not converge (impossibility to reach the discrepancy principle criterion). When it converges, the THS rebuilt by Branch modes remains far from the exact one. In the case of Fig.11-c (s.t.d of the noise $\sigma = 0.02$), the optimal number of branch modes was found equal to $N_m = 29$ and removing 11 high frequency modes still produces an oscillating reconstruction while not modifying really the temperature residuals. In Fig.11-e (s.t.d of the noise $\sigma = 0.1$), the relative errors are of the order of 10% in the linear parts of source and knowing that the optimal number of branch modes found in that case is precisely 18 (equal to the number of optimal Fourier modes in the same case), the comparison with Fig.11-f definitely shows that the Branch modes is too sensitive to noise on the observation.

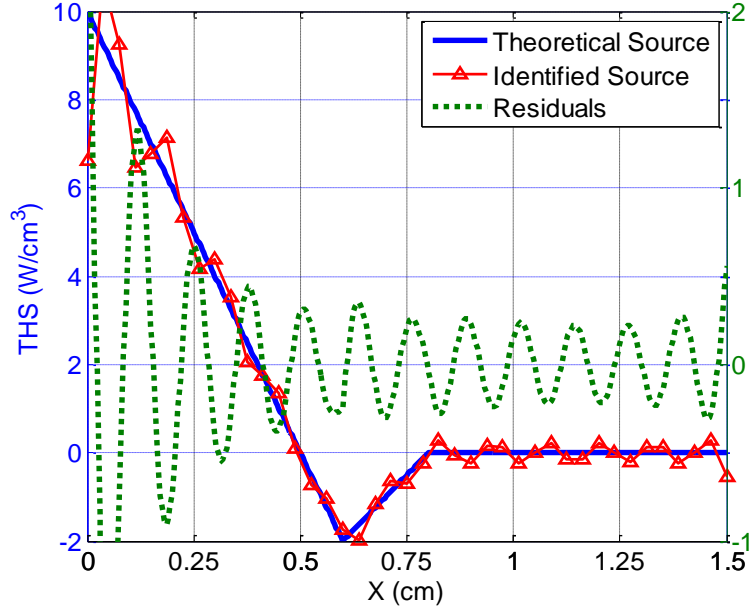

Fig 12. Test case 1 : THS reconstructed at $t = 20s$ through 25 Branch modes without noise.

This trend is confirmed –in the presence of a small amount of noise- when the THS reconstruction is based on the few first identified experimental modes. For this and using the same TSVD or KLD scheme as the one used by operating Eqs (29) to (36), but applied here to the Gram matrix of the experimental signal (Barrio 2012), it is possible to calculate the empirical SVD modes of the signal. The temperature field (or profiles) can be decomposed on these singular modes. It was found by simulations that only the same few first modes concentrate most of the energy of the signal (>95% for only 3 modes!): these modes remain un-affected by the noise of the signal up to a given number. The coefficients of the source decomposition on this set of singular modes can be straightforwardly obtained without any extra regularization from an explicit scheme applied to the differential system of the type of eq.(5). Such a reconstruction is shown on Fig. 13-a. Based on this initialization of the THS (first three component of the states vector $\boldsymbol{B}^{(0)}$), the branch modes works now better (Fig. 13-b) in the case of non negligible superimposed noise. Although it still not compares as excellent as Fourier's modes (compare to Fig. 11-c), it produces reasonably good estimations. In conclusion, the Fourier basis appears to have better performances in this IHSP. It reconciles a good spreading of state values in the frequency domain (ensuring stability of the inversion process) with a sufficiently great number of modes to ensure good reconstructions of locally sharp functions.

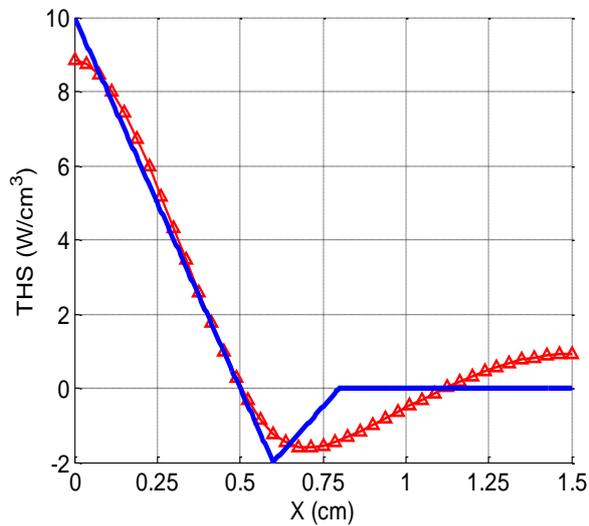 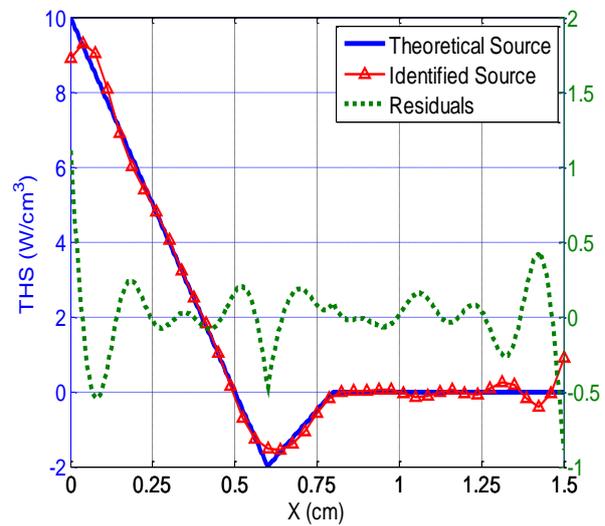

Fig 13-a. A priori information of THS at $t = 20s$ estimated from 3 empirical singular modes - $\sigma = 0.02$

Fig. 13-b. Reconstructed THS at $t = 20s$ from initial THS profile of Fig13-a - 18 Branch modes - $\sigma = 0.02$

## 6  Conclusion

It is well known that model reduction provides an efficient numerical method for solving inverse heat conduction problems. This was done here on a particular Inverse Heat Source Problem using the Branch basis. This is a more general basis of decomposition which is very performing for the reduction of very complex heat transfer problems. It has been investigated here to see if it could enhance the general performances of a given inverse algorithm compared to classical Fourier basis. In that case of a "simple" 1-D heat transfer problem, the answer is no and it is one strong conclusion here to state that, because we have a diffusion and (non constant) advection problem, the best basis for reconstruction cannot be found as the solution of the exact linear eigenvalue problem. As a result, it has to be chosen in some empirical manner with respect to its performances regarding the considered inverse problem (expected form of source terms or boundary conditions of the problem). The second main issue of this work lies in the regularizing treatment brought to the "adjoint-conjugate gradient" optimization. Within the framework of the CGM, it has been shown that KLD or TSVD applied to the conjugate directions matrix has a strong regularizing effect which can be advantageously used for any problem.

# Appendix

## A. Branch Basis Decomposition

In this section, we detail the Branch basis approach that was used to decompose both the temperature and the heat source. This basis is determined by the following generalized eigenvalue problem (A.1) where the specificity of the Branch basis lies in the generalized Steklov boundary condition (A.2).

$$\forall X \in \bar{\Omega} = [0, L] \qquad k \frac{\partial^2 V_i(X)}{\partial X^2} - cv_0 \frac{\partial V_i(X)}{\partial X} = \lambda_i c V_i(X) \qquad (A.1)$$

$$-k \frac{\partial V_i(X)}{\partial X}\bigg|_{X=0} + cv_0 V_i(X)|_{X=0} = -\lambda_i \zeta V_i \text{ and } -k \frac{\partial V_i(X)}{\partial X}\bigg|_{X=L} = \lambda_i \zeta V_i \qquad (A.2)$$

$\{V_i, \lambda_i\}$ stands for the $i^{th}$ eigenvector and the associated $i^{th}$ eigenvalue. The parameter $\zeta$ which appears in (A.2) is constitutive of Steklov boundary conditions and ensures dimensional homogeneity of eigenvalues $\lambda_i$. In 1D, it can be demonstrated that this coefficient must be chosen as $\zeta = c * L/2$. When the advection term is taken into account for the Branch basis definition, this eigenvalue problem is no longer self-adjoint. The spectral Branch basis can be determined when:

**i.** a single "entering" advective flux is considered in the Steklov boundary equations

**ii.** both proper and adjoint eigenmodes are considered (condition to form a bi-orthogonal basis). For the adjoint eigenvalue equation, only an output advective flux is considered.

The adjoint problem (star superscript) is now given by

$$\forall X \in \bar{\Omega} = [0, L] \qquad k \frac{\partial^2 V_i^*}{\partial X^2} + cv_0 \frac{\partial V_i^*}{\partial X} = c\lambda^* V_i^* \qquad (A.3)$$

$$-k \frac{\partial V_i^*}{\partial X}\bigg|_{X=0} = -\lambda^* \zeta V_i^* \quad and \quad -k \frac{\partial V_i^*}{\partial X}\bigg|_{X=L} = (\lambda^* \zeta + cv_0) V_i^* \qquad (A.4)$$

The set of functions $(V_i, V_j^*)$ do form a basis, respecting the bi-orthogonality condition (Neveu, 2005):

$$C_D(V_i, V_j^*) = \int_D cV_i\overline{V_j^*}dX + \zeta(V_i(0)\overline{V_j^*}(0) + V_i(L)\overline{V_j^*}(L)) = \delta_{ij} \tag{A.5}$$

For differential one-directional eigen problems (A.1) and (A.3), proper and adjoint eigenmodes can be computed semi-analytically. In dimensionless form and dropping the i index for clarity, (A.1) turns into

$$\frac{\partial^2 V}{\partial x^2} - Pe\frac{\partial V}{\partial x} = \lambda V \tag{A.6}$$

$$-\frac{\partial V}{\partial x}\bigg|_{x=0} = -(\lambda\zeta + Pe)V \quad \text{and} \quad -\frac{\partial V}{\partial x}\bigg|_{x=1} = \lambda\zeta V \tag{A.7}$$

where $x = \frac{X}{L}$, $Pe = \frac{cv_0 L}{k}$, $\lambda = \frac{cL^2}{k}\lambda$, $\zeta = \frac{1}{2}$ and $x \in [0,1]$. Looking for a solution in the form $V(x) = e^{rx}$, the associated characteristic equation $r^2 - Pe.r - \lambda = 0$ leads to the roots

$$r_1 = Pe + \frac{\sqrt{Pe^2 + 4\lambda}}{2}, \quad r_2 = Pe - \frac{\sqrt{Pe^2 + 4\lambda}}{2}$$

By setting, $\alpha = \frac{Pe}{2}$, $\beta = \frac{\sqrt{Pe^2 + 4\lambda}}{2}$ and thus $\lambda = \beta^2 - \frac{Pe^2}{4}$, the solution of the problem (A.6) is:

$$V(x) = e^{\alpha x}(A\cosh(\beta x) + B\sinh(\beta x)) \tag{A.8}$$

The coefficients will be determined through the boundary conditions (A.7), which give:

$$(\lambda + \alpha)A = \beta B \tag{A.9}$$

$$\left(\alpha + \frac{\lambda}{2} + \beta\tanh(\beta)\right)A = -\left((\alpha + \frac{\lambda}{2})\tanh(\beta) + \beta\right)B \tag{A.10}$$

As a result, we have a transcendental equation related to $\beta$:

$$\tanh(\beta) = \frac{(Pe + \lambda)\beta}{-\beta^2 - (\lambda + Pe)^2/4} \tag{A.11}$$

For $\beta$ being real, we have at most two eigenfunctions defined by (A.8), whereas an infinite and countable set of eigenfunctions is obtained for the imaginary $\beta_i = jq_i, \forall q_i \in \mathbb{R}$:

$$V_i(x) = e^{\alpha x}\left(\cos(q_i x) + \frac{\left(-q_i^2 - \alpha^2\right)/2 + \alpha}{q_i}\sin(q_i x)\right) \tag{A.12}$$

Respectly, the adjoint eigenfunctions are given as

$$V_i^*(x) = e^{-\alpha x}\left(\cos(q_i x) + \frac{\left(-q_i^2 - \alpha^2\right)/2 + \alpha}{q_i}\sin(q_i x)\right) \tag{A.13}$$

One can easily find out that $V_i$ and $V_i^*$ are identical except the exponential function at left. They are real-valued but the orthogonality condition is now verified to be of the order of $10^{-5}$ (with a mesh grid of 400 nodes). Note that if $v_0$ is set equal to 0, one recovers the pure diffusion branch eigenmodes (self-adjoint operator).

As a result, the temperature can be written as

$$T(X,t) = \sum_{i=1}^{N} z_i(t) V_i(X) \tag{A.14}$$

And, in a similar way, the source can be written as

$$q(X,t) = \sum_{i=1}^{N} b_i(t) c V_i(X). \tag{A.15}$$

The inverse problem related to eq (5) is now formatted so as to produce an estimation of coefficients $z_i(t)$ and $b_i(t)$ which are related through the following system of state ODEs:

$$\forall m \in 1, \dots, N$$

$$\sum_{i=1}^{N} \dot{z}_i(t)\gamma_{im} = \lambda_m z_m(t) - \sum_{i=1}^{N} z_i(t) P_{im}(t) + Q_m(t) + \sum_{i=1}^{N} \dot{b}_i(t)\gamma_{im} \tag{A.16}$$

with

$$\gamma_{im} = \delta_{im} - \zeta(V_i(L)\overline{V_m^*}(L) + V_i(0)\overline{V_m^*}(0))$$

$$Q_m(t) = -\varphi_2(t)\overline{V_m^*}(L) + \varphi_1(t)\overline{V_m^*}(0)$$

$$P_{im}(t) = \int_D cv(X,t)\frac{\partial V_i}{\partial X}\overline{V_m^*}dX - cv_0 V_i(0)\overline{V_m^*}(0)dX$$

## B: Fourier Basis Decomposition:

This approach is identical but allows for a different treatment of the Boundary Condition.

We first consider that the temperature variable solving problem (1) can be separated in a homogeneous component $T_h$ and a particular function accounting for the BC,

$$T(X,t) = -\frac{2L}{k\pi}\varphi_1(t)\sin\left(\frac{\pi X}{2L}\right) + \frac{2L}{k\pi}\varphi_2(t)\cos\left(\frac{\pi X}{2L}\right) + T_h(X,t) \tag{B.1}$$

From (B.1), we have

$$\frac{\partial T}{\partial t} = -\frac{2L}{k\pi}\dot{\varphi}_1(t)\sin\left(\frac{\pi X}{2L}\right) + \frac{2L}{k\pi}\dot{\varphi}_2(t)\cos\left(\frac{\pi X}{2L}\right) + \dot{T}_h(X,t) \tag{B.2}$$

$$\frac{\partial T}{\partial X} = -\frac{1}{k}\varphi_1(t)\cos\left(\frac{\pi X}{2L}\right) - \frac{1}{k}\varphi_2(t)\sin\left(\frac{\pi X}{2L}\right) + \frac{\partial T_h}{\partial X} \tag{B.3}$$

$$\frac{\partial^2 T}{\partial X^2} = \frac{\pi}{2kL}\varphi_1(t)\sin\left(\frac{\pi X}{2L}\right) - \frac{\pi}{2kL}\varphi_2(t)\cos\left(\frac{\pi X}{2L}\right) + \frac{\partial^2 T_h}{\partial X^2} \tag{B.4}$$

Replacing (B.2,3,4) in the system (1) lead to a homogeneous form in terms of $T_h$ variable:

$$c\frac{\partial T_h}{\partial t}(X,t) = k\frac{\partial^2 T_h}{\partial X^2}(X,t) - cv(X,t)\frac{\partial T_h}{\partial X}(X,t) + p(X,t) + q(X,t) \tag{B.5}$$

$$-k\frac{\partial T_h}{\partial X}\bigg|_{X=0} = 0 \tag{B.6}$$

$$-k\frac{\partial T_h}{\partial X}\bigg|_{X=L} = 0 \tag{B.7}$$

$$T_h(X,0) = T_0 + \frac{2L}{k\pi}\varphi_1(0)\sin\left(\frac{\pi X}{2L}\right) - \frac{2L}{k\pi}\varphi_2(0)\cos\left(\frac{\pi X}{2L}\right) \tag{B.8}$$

with

$$p(X,t) = \frac{2cL}{k\pi}\left(\dot{\varphi}_1(t)\sin\left(\frac{\pi X}{2L}\right) - \dot{\varphi}_2(t)\cos\left(\frac{\pi X}{2L}\right)\right)$$

$$+ \frac{\pi}{2L}\left(\varphi_1(t)\sin\left(\frac{\pi X}{2L}\right) - \varphi_2(t)\cos\left(\frac{\pi X}{2L}\right)\right)$$

$$+ \frac{cv(X,t)}{k}\left(\varphi_1(t)\cos\left(\frac{\pi X}{2L}\right) + \varphi_2(t)\sin\left(\frac{\pi X}{2L}\right)\right)$$

With Neuman conditions, the function $T_h$ defined in $[0,L]$ can be decomposed by a cosine set, such that

$$T_h = \frac{z_0(t)}{2} + \sum_{i=1}^{\infty} z_i(t)\cos\left(\frac{i\pi X}{L}\right), \tag{B.9}$$

as well as the source term

$$q(X,t) = \frac{b_0(t)}{2} + \sum_{i=1}^{\infty} b_i(t)\cos\left(\frac{i\pi X}{L}\right) \tag{B.10}$$

Thanks to the natural orthogonal property of cosine modes, the same mathematical problem as in eqs (B.9,10) can be settled up in terms of states $z_i(t)$ and $b_i(t)$. The drawback of this approach is that it requires to calculate the derivatives of the fluxes $\dot{\varphi}_1(t), \dot{\varphi}_2(t)$, which is a risk of increasing the instability of the problem because these fluxes will not be perfectly known.

## C. Algebraic discretized form of state equations

Given a set of basis functions $\{V_i\}_{i \in \mathbb{N}}$ -and independently of the selected spectral basis- the system (1) can be turned into (5,6) through simple scalar product in $L^2(\Omega)$ space, which results in the following matrix coefficients

$$A(t) = (\langle k \frac{\partial^2 V_i}{\partial X^2}, V_j \rangle_\Omega - \langle cv(X,t) \frac{\partial V_i}{\partial X}, V_j \rangle_\Omega)_{ji}$$

$$C = c\mathbf{I}$$

$$D = (\langle V_i, V_j \rangle_\Omega)$$

Actually, the forms of $A(t), C, D$ can change since the orthogonal property of the basis could be defined differently. Note that if $\langle V_i, V_j \rangle_\Omega = \delta_{ij}$, we have $D = I$.

In the case of the Fourier basis decomposition of annex B, we have

$$A = diag\left(0, -k\left(\frac{\pi}{L}\right)^2, -k\left(\frac{2\pi}{L}\right)^2, \dots, -k\left(\frac{(N_m-1)\pi}{L}\right)^2\right) - P$$

$$P = \begin{pmatrix} 0 & \cdots & \frac{2(N_m-1)\pi}{L^2} \int_0^L cv(X,t) \sin\left(\frac{(N_m-1)\pi X}{L}\right) dX \\ \vdots & & \vdots \\ 0 & \cdots & \frac{2(N_m-1)\pi}{L^2} \int_0^L cv(X,t) \sin\left(\frac{(N_m-1)\pi X}{L}\right) \cos\left(\frac{(N_m-1)\pi X}{L}\right) dX \end{pmatrix}$$

$$C = diag(c, c, \dots, c)$$

$$D = I$$

# References

V. Isakov, *Inverse Source Problems,* American Mathematical Society, Providence, (1990).

Auffray, N.; Bonnet, M.; Pagano, S., Identification of transient heat sources using the reciprocity gap, Inverse Problems in Science and Engineering, 21(4), 721-738, (2013).

A. Chrysochoos and F. Belmahjoub, Thermographic analysis of thermo-mechanical couplings, Arch. Mech., 44, 55–68 (1992).

S. André, N. Renault, Y. Meshaka, C. Cunat, From the thermodynamics of constitutive laws to thermomechanical experimental characterization of materials: An assessment based on inversion of thermal images, Continuum Mechanics and Thermodynamics: 24(1), 1-20, (2012), doi:10.1007/s00161-011-0205-xA.

Vincent, L., On the ability of some cyclic plasticity models to predict the evolution of stored energy in a type 304L stainless steel submitted to high cycle fatigue, European Journal of Mechanics A/Solids, 27, 161–180, (2008)

C. Doudard, S. Calloch, F. Hild and S. Roux, Identification of heat source from infrared thermography: Determination of 'self-heating' in a dual-phase steel by using a dog bone sample, Mechanics of Materials, 42, 55-62, (2010).

Boulanger, T., Chrysochoos, A., Mabru, C., Galtier, A., 2004. Calorimetric analysis of dissipative and thermoelastic effects associated with the fatigue behavior of steels. Int. J. Fatigue 26, 221–229.

N. Renault, S. André, D. Maillet and C. Cunat, A two-step regularized inverse solution for 2-D heat source reconstruction, International Journal of Thermal Sciences 47 (2008) 834-847.

N. Renault, S. André, D. Maillet and C. Cunat, A spectral method for the estimation of a thermomechanical heat source from infrared temperature measurements, International Journal of Thermal Sciences, 49 (2010), 1394-1406.

Rouquette, S, Guo, J, Le Masson, P, Estimation of the parameters of a Gaussian heat source by the Levenberg-Marquardt method: Application to the electron beam welding, Int. J Thermal Sciences, 46(2), 128-138, (2007). DOI: 10.1016/j.ijthermalsci.2006.04.015.

Yi, Z Murio, DA, Source term identification in 1-D IHCP, Computers & Mathematics with applications, 47(12), 1921-1933, (2004). DOI: 10.1016/j.camwa.2002.11.025.

Delpueyo, D., Balandraud, X., Grédiac, M., Heat source reconstruction from noisy temperature fields using an optimised derivative Gaussian filter, Infrared Physics & Technology, 60, 312–322, (2013).

Farcas, A., Lesnic, D., The boundary-element method for the determination of a heat source dependent on one variable, Journal of Engineering Mathematics, 54(4), 375-388, (2006).

Le Niliot, C, Rigollet, F, Petit, D, An experimental identification of line heat sources in a diffusive system using the boundary element method, International Journal of Heat and Mass Transfer, 43(12), 2205-2220, (2000). DOI: 10.1016/S0017-9310(99)00285-9 .

Liu, CS, An iterative algorithm for identifying heat source by using a DQ and a Lie-group method, Inverse Problems in Science and Engineering, 23(1), 67-92 (2015). DOI: 10.1080/17415977.2014.880907.

Hasanov, A., Pektas, B., Identification of an unknown time-dependent heat source term from overspecified Dirichlet boundary data by conjugate gradient method, Computers & Mathematics with Applications, 65(1), 42-57, (2013). DOI: 10.1016/j.camwa.2012.10.009 .

Erdem, A.; Lesnic, D.; Hasanov, A., Identification of a spacewise dependent heat source, Applied Mathematical Modelling, 37(24), 10231-10244, (2013).

Maalej, T, Maillet, D., Fontaine, J-R, Estimation of position and intensity of a pollutant source in channel flow using transmittance functions, Inverse Problems, 28(5), 055010, (2010). doi:10.1088/0266-5611/28/5/055010

Rap, A, Elliott, L., Ingham, D.B., Lesnic, D., Wen, X., The inverse source problem for the variable coefficients convection-diffusion equation, Inverse Problems in Science and Engineering, 15(5), 413-440, (2007). DOI: 10.1080/17415970600731274

De Sousa, D.M., Roberty, N.C., An inverse source problem for the stationary diffusion–advection–decay equation, *Inverse Problems in Science and Engineering*, 20(7), 891–915, 2012.

Videcoq, E., Quemener, O., Lazard, M., Neveu, A., Heat source identification and on-line temperature control by a Branch Eigenmodes Reduced Model, International Journal of Heat and Mass Transfer, 51, 4743–4752 (2008).

Loeven

Hanson, F.B., *Applied Stochastic Processes and Control for Jump-Diffusions: Modeling, Analysis and Computation*, Appendix A, Society for Industrial and Applied Mathematics, 2007.

Liberzon, D., *Calculus of variations and optimal control theory, UIUC Courses,* (2011).

Jonathan Richard Shewchuk, An Introduction to the Conjugate Gradient Method Without the Agonizing Pain, (1994).

Hanke-Bourgeois, M., *Conjugate gradient type methods for ill-posed problems,* CRC Press, (1995).

Shenefelt, J.R, Luck, R., Taylor, R.P., Berry, J.T., Solution to inverse heat conduction problems employing singular value decomposition and model reduction, International Journal of Heat and Mass Transfer, 45(1), 67-74, (2002).

Park, H.M., Jung, W.S., The Karhunen-Loève Galerkin method for the inverse natural convection problems, International Journal of Heat and Mass Transfer, 44,155-167, (2001).

Pujol J., The solution of nonlinear inverse problems and the Levenberg-Marquardt method, *Geophysics* (SEG), 72(4), W1–W16, (2007). doi:10.1190/1.2732552.

Barrio, E.P.D., Dauvergne, J.L., Pradere, C. Thermal characterization of materials using Karhunen-Loève decomposition techniques, Inverse Problems in Science & Engineering, 20(8), 1115-1143, (2012).

Cabeza, J.M.G, Garcia, J.A.M , Rodriguez, A.C. , A sequential algorithm of inverse heat conduction problems using singular value decomposition, International Journal of Thermal Sciences, 44, 235-244, (2005).

Ye, J. Andre, S., Farge, L., Kinematic Study of Necking in a SemiCrystalline Polymer through 3D Digital Image correlation, Int. J. Solids & Structures, dx.doi.org/10.1016/j.ijsolstr.2015.01.009, published on-line 02-02-2015.

Versteeg, H.K., Malalasekera, W., *Introduction to Computational Fluid Dynamics: The Finite Volume Method,* PEARSON Prentice Hall, $2^{nd}$ (2007).

Hochstenbach, M.E., Reichel, L., An Iterative Method for Tikhonov Regularization with A General Linear Regularization Operator, Journal of Integral Equations and Applications, 22(3), 465-482, (2010).

Neveu, A., Approche modale pour les systèmes thermiques, Internal Note - Laboratoire de Mécanique et d'Energétique d'Evry, *Personal Communication,* 2005.